# Analytic solution to swing equations in power grids with ZIP Load Models

# Analytic solution to swing equations in power grids with ZIP load models

HyungSeon Oh[1*]

[1] Department of Electrical and Computer Engineering, United States Naval Academy, Annapolis, Maryland, USA
* Email: hoh@usna.edu

Abstract

**Objective**: This research pioneers a novel approach to obtain a closed-form analytic solution to the nonlinear second order differential swing equation that models power system dynamics. The distinctive element of this study is the integration of a generalized load model known as a *ZIP* load model (constant impedance *Z*, constant current *I*, and constant power *P* loads).

**Methods**: Building on previous work where an analytic solution for the swing equation was derived in a linear system with limited load types, this study introduces two fundamental novelties: 1) the innovative examination and modeling of the *ZIP* load model, successfully integrating constant current loads to augment constant impedance and constant power loads; 2) the unique derivation of voltage variables in relation to rotor angles employing the holomorphic embedding (*HE*) method and the Padé approximation. These innovations are incorporated into the swing equations to achieve an unprecedented analytical solution, thereby enhancing system dynamics. Simulations on a model system were performed to evaluate transient stability.

**Results**: The *ZIP* load model is ingeniously utilized to generate a linear model. A comparison of the developed load model and analytical solution with those obtained through time-domain simulation demonstrated the remarkable precision and efficacy of the proposed model across a range of IEEE model systems.

**Conclusion**: The study addresses the key challenges in power system dynamics, namely the diverse load characteristics and the time-consuming nature of time-domain simulation. Breaking new ground, this research proposes an analytical solution to the swing equation using a comprehensive *ZIP* model, without resorting to unphysical assumptions. The close-form solution not only assures computational efficiency but also preserves accuracy. This solution effectively estimates system dynamics following a disturbance, representing a significant advancement in the field.

I. NOMENCLATURE

*Nb*    number of buses

$E_{ref}$    reference voltage magnitude

$E_I$    internal voltage magnitude vector at *Ibus*

*Ibus*    transformed internal generator nodes

IEEE    Institute of Electrical and Electronics Engineers

*Kbus*    set of nodes directly connected to a generator

*Mbus*    set of nodes that are not directly connected to a generator

*KM*    union of *Kbus* and *Mbus*, $K \cup M$

*NI*    number of nodes in *Ibus*

*NKM*    number of nodes in *Kbus* or in *Mbus*

$W_j$    reciprocal voltages at the $j^{th}$ node, $W_j = 1/v_j$

$Y_{bus}$    nodal admittance matrix

$Y_{tr}$    admittance matrix where $Y_{tr} = G_{tr} + \mathbf{j}B_{tr}$

**j**    $\sqrt{-1}$

$v_x$    real component of voltage vector

$v_y$    imaginary component of voltage vector

$v$    complex voltage vector, $v = v_x + \mathbf{j}v_y$

\*    complex conjugate

## II. INTRODUCTION

As societal demand for electric power grows, maintaining synchronization among synchronous generators becomes increasingly difficult. Failure to maintain synchronization might lead to fluctuations in the rotational motion of the generators over time. The swing equations [1], which are complex nonlinear second-order differential equations, govern the behavior of the synchronous generators. The problem has remained unanswered for almost a century, and research on rotor motions post-disturbance has primarily concentrated on either physically unacceptable assumptions [2]–[9] or time-consuming numerical simulations [10]–[16], which may not be suited for assessing real-time system stability.

An analytical solution for linear loads such as frequency-dependent loads and constant impedance loads was recently developed, allowing power grid stability to be assessed with minimal computational effort [17]. However, this approach does not take into account nonlinear loads, such as constant power loads (PQ loads). The BIG model has been proposed [18][19], but it implies load linearization and a minimal deviations from the static operating point, which may not be valid during disturbances. The deviation could lead to substantial errors in the simulation results [10][11][16]. Please see the surveys in [20], [21] for further information on load modeling.

During the transients, voltages and flows are determined by power flow equations. To determine the voltages for a given set of *PQ*, *PV*, and slack buses, a power flow (*PF*) problem is solved. Because the *PF* problem is *NP*-hard (non-deterministic polynomial-time), finding a precise solution within an acceptable timeframe is challenging. To solve the *PF* problem, heuristic approaches, such as Newton-Raphson or Gauss-Seidel methods [1], are typically utilized. Each node usually has two known values and two equations for calculating the real and imaginary parts of voltages, as well as real and reactive power injections. However, a novel tensor computation-based method has been presented that does not require uniform knowledge and assigns *2N* knowns from all nodes [22].

Holomorphic embedding (*HE*) is another method being investigated for an analytical solution to the *PF* problem [23][24]. *HE* entails equating the *s* term by term and extending the voltages, reciprocal voltages, and reactive power generation in terms of the variable *s*. The zeroth order solution, or "germ solution," is determined independently of the loading conditions, although the first and higher order terms are affected by them. The *HE* method is extended to numerical dynamic simulations, albeit with the strict requirements that the variables be approximated as a truncated sum of McLaurin series [25]. Because this assumption is only valid within an unknown convergence radius, the number of terms to include after truncation is determined heuristically, which potentially leads to inaccurate results if the number of terms is incorrectly assigned.

In this study, the HE method will be employed to derive voltage interdependence for various load characteristics. The following is how the paper is structured: Section III provides a load modeling approach for integrating load characteristics into swing equations in an efficient and precise manner. Section IV shows how nodal voltages are linearly proportional to to rotor angles (i.e., voltages at *Ibus*). Section V provides the analytical solution to a differential equation established inside a validity region in order to approximate swing equations using nonlinear load models. Section VI extends the analytical solution's application beyond its validity region. Section VII discusses the numerical methodologies that were employed in this study. Section VIII describes the numerical results obtained by simulating dynamics with *ZIP* loads and explains stability assessments for selected examples. Finally, Section IX summarizes the study's findings and draw conclusions.

## III.  LOAD MODELING

Because the solution to swing equation stated in [17] already accounts for frequency dependent loads, this section will focus on a frequency-independent load model. The *ZIP* model, which incorporates constant impedance *Z*, constant current *I*, and constant power *P* loads, is a popular load model for power flow and dynamic simulations [20], [21]. This load model accurately represents load behavior under normal and fault conditions, making it a useful tool for power flow studies and stability analysis. Constant impedance *Z*, constant current *I*, and constant power *P* loads are addressed in detail in the sections that follow.

*III.1. Constant Impedance Loads*

Constant impedance (*Z*) loads have steady impedance characteristics, which means their electrical resistance remains constant regardless of the applied voltage. As a result, the current delivered to the load is directly proportional to the voltage, however the power delivered to the load increases with the square of the magnitude of the voltage. Constant impedance loads accurately depict resistance-based loads such as electric heating devices and incandescent light bulbs.

One approach to model constant impedance loads is to embed them in the $Y_{bus}$ matrix. The $j^{th}$ element in the diagonal $Y_{bus}^{add}$ matrix in this method is $Y_{bus}^{add}(j,j) = \frac{S_j^{0*}}{|v_j|^2}$ [26], where $S_j^{0*}$ represents the complex conjugate of the constant impedance loads at the $j^{th}$ node when the voltage magnitude at the node $|v_j|$ is 1.0 per unit (p.u.). When the voltage magnitude is $|v_j|$, the actual realized power injection at Node $j$ is $S_j^* = (v_j i_j^*)^* = Y_{bus}^{add}(j,j)|v_j|^2 = S_j^{0*}|v_j|^2$ [26]. Consequently, the modified nodal admittance matrix is $Y_{bus} = Y_{bus}^0 + diag(S_j^{0*})$, where $Y_{bus}^0$ is the nodal admittance matrix without Z loads. The modified matrix $Y_{bus}$ is then used to calculate the network current $I_{bus}$ with constant impedance, which can be estimated using the formula $I_{bus} = Y_{bus} v$.

*III.2. Constant Current Loads*

Constant current loads are electrical loads that provide power proportional to the magnitude of the voltage at the node, $S_j = S_j^0 |v_j|$ [26]. It is, however, difficult to establish the relationship between current and voltage based on the linear relationship between power and voltage magnitude. As a result, for modeling purposes, we approximated constant current loads in terms of other load types. During transients, the voltage magnitudes $|v_j|$ at the KM nodes do not deviate significantly from known values and may be expressed as $|v_j| = |v_j^0| + \Delta_j$, where $|\Delta_j| \ll |v_j^0|$. Taylor's series expansion approximates $|v_j|^2$, $|v_j|^2 \cong |v_j^0|^2 + 2|v_j^0|\Delta_j$, allowing constant current loads to be decomposed into constant impedance and constant power loads, $S_j = S_j^0|v_j| = \left(\frac{S_j^0}{2|v_j^0|}\right)|v_j|^2 + \left(\frac{S_j^0}{2}\right)|v_j^0|$. The deviation of the approximation $S_j^{app}/S_j^{true}$ is calculated as an error measurement; $S_j^{app}/S_j^{true} = \frac{1}{2}\frac{|v_j|}{|v_j^0|} + \frac{1}{2}\frac{|v_j^0|}{|v_j|}$.

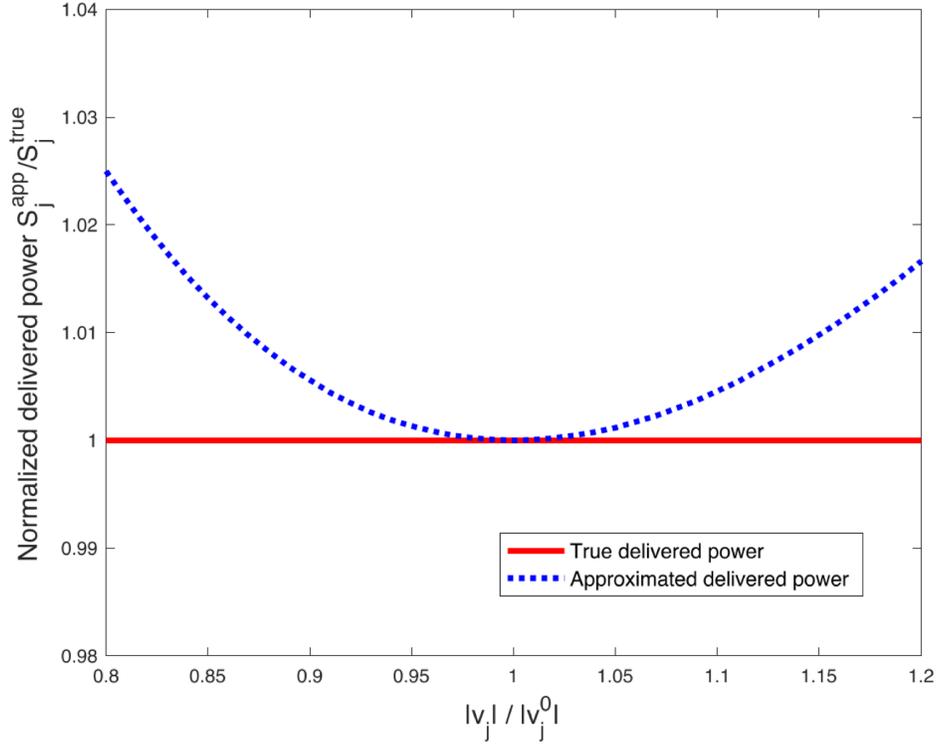

Figure 1. The normalized delivered power $S_j^{app}/S_j^{true}$ in terms of the normalized voltage magnitude $|v_j|/|v_j^0|$. The deviations of the approximation are 2-3% within 20% deviation of the voltage magnitudes.

Figure 1 depicts the true delivered normalized power to constant current loads and its approximation, which contains two quadratic components in relation to voltage magnitude and do not vary with magnitude. Figure 1 shows that the approximation closely matches the true values within a range of 2-3% for a wide range of voltage magnitudes, with the approximation deviating by around 20% from the voltage magnitude at which it was made. If the voltage magnitudes exceed a specific range, the approximation may no longer be valid, and an updated approximation with the new voltage magnitude $|v_j^0|$ may be necessary. However, there was little need for this in our stability studies.

*III.3. Constant Power Loads*

A constant power load is an electrical load that requires a consistent power supply to operate regardless of voltage or current fluctuations. Although the load's impedance varies, the power consumed remains constant. Electric motors, fluorescent lighting, and induction heating are examples of such loads. Voltage-dependent loads, on the other hand, are easier to support on the transmission network since they respond to variations in system voltage. Constant impedance loads, for example, fluctuate proportionally to the square of the voltage, resulting in a greater reduction when the voltage falls below 1.0 p.u.. This constant power load type is observed in numerous power flow analyses, such as MATPOWER [27] and Holomorphic Embedding [23].

*III.4. Concluding Remark on the ZIP Load Models*

Throughout this study, the *ZIP* load models are used in a variety of ways. Constant current *I* loads are first transformed into constant impedance and constant power loads. The extended *Z* and *P* loads are created by combining these transformed components with the existing *Z* and *P* loads. The extended constant impedance loads *Z* are incorporated into the nodal admittance matrix $Y_{bus}$, modifying the diagonal components of the nodal admittance matrix $Y_{bus}^0$, which is generated using the branch data and shunt elements. *P* continues to be a load for the power flow problem or quasi-power flow problem (*QPF* as defined in Section IV), which is discussed in the following section. As a result, the *PF* problem, or *QPF*, establishes the voltages required to fulfill the extended *P* loads throughout the system, using the modified $Y_{bus}$.

IV. LINEARIZED INTERDEPENDENCE BETWEEN VOLTAGES

Real and reactive power injections, as well as real and imaginary voltages, are required to define the state of a node. Traditional *PF* problems use two known variables at each node to determine the other two, which vary depending on the bus type: 1) slack bus with known voltages; 2) *PV* bus with specified real power injection and voltage magnitude; and 3) *PQ* bus with both real and reactive power injections specified. The power balance equations at Bus *i*, denoted as $e_i^T Y_{bus} v = \frac{1}{v_i^*}(p_i - \mathbf{j}q_i)$ [1], and the computation of the other two variables using the two known.

During transients, the voltage magnitudes at *Ibus* stay constant [28]. The goal of *QPF* is to determine the voltages at *KM* in response to the rotor angle (equivalent to the voltage angle at *Ibus*). This research addresses the *QPF* problem, which is a form of *PF* problem in which the remaining two variables are calculated using power balance equations and two known variables. The *QPF* formulation differs slightly from standard *PF* problems in that it has two types of buses: reference buses with known real and imaginary voltages (all *Ibus*) and *PQ* buses.

Load characteristics are critical and this study takes *ZIP* loads into account. After the decomposition of *I* loads, *QPF* determines voltages throughout the system that satisfy the extended *P* loads with the modified nodal admittance matrix in terms of the extended *Z* loads. The Kronecker product is capable of efficiently solving the *QPF* problem [22]. Kirchhoff's current law finds $I_j = e_j^T(Y_{bus} v + I_j^0) = \frac{S_j^*}{\|v_j\|^2} e_j^T v$ at Bus $j \in$ *KM* [1]. In this case, $I_j^0$ is a designated current for all *j* in *KM*, serving as a placeholder to improve load model flexibility, such as shunt currents, creating a linear relationship between the voltages. Previous investigations have shown that when only *Z* loads are present, the linear relationship is accurate [17]. Many attempts have been made to derive a linear relationship between nodal voltages. The linear relationship has the advantage of allowing the power flow equations to be concisely expressed in a smaller number of nodes, which is why it is sometimes referred as a network reduction. Section IV.1. illustrates conventional techniques, and Section IV.2. presents the proposed approach based on Holomorphic Embedding.

*IV.1. Conventional Methods to Represent $v_{KM}$ in terms of $v_I$*

The Ward reduction method [29] is employed when the voltage magnitudes at *KM* are almost constant,. The current at *KM* can be expressed as a linear combination of the voltages at *KM* and *Ibus*, $I_{KM} = Y_{bus}^{KM} v_{KM} + Y_{bus}^{I} v_I = diag(S_{KM}^*/|v_{KM}|^2) v_{KM} + I_{KM}^0$ where $Y_{bus}^{KM}$ is a partitioned $Y_{bus}$ matrix with the rows and columns correspond to *KM*; $Y_{bus}^{I}$ is a partitioned $Y_{bus}$ matrix in which the rows correspond to *KM* but the columns correspond to *Ibus*; $S_{KM}^*/|v_{KM}|^2$ is the complex conjugate of the delivered power to *KM* normalized with the squared voltage magnitude; and $I_{KM}^0$ is the designated current at *KM*.

The following equation can estimate the voltages at *KM* using the known voltages at *Ibus* and the power injections, assuming that the voltage magnitude at *KM* remains constant: $v_{KM} = \{-[Y_{bus}^{KM} - diag(S_{KM}^*/|v_{KM}|^2)]^{-1} Y_{bus}^{I}\} v_I + \{[Y_{bus}^{KM} - diag(S_{KM}^*/|v_{KM}|^2)]^{-1} I_{KM}^0\}$ [29]. This method can potentially reduce the complexity of the power flow problem while increasing computational efficiency.

Voltage sensitivities close to the original dispatch can also be investigated: $f(v) = f(v_{KM}, v_I) \cong f(v^0) + [\nabla_{v_{KM}} f]^T \Delta v_{KM} + [\nabla_{v_I} f]^T \Delta v_I$. The nodal power balance, $f(v) = f(v^0) = 0$, must be maintained while voltages fluctuate. As a result, the first-term Taylor series expansion provides a linear relationship between *Ibus* and other buses' voltages: $v_{KM} = \left(-[\nabla_{v_{KM}} f]^{-T} [\nabla_{v_I} f]^T\right) v_I + \left(v_{KM}^0 + [\nabla_{v_{KM}} f]^{-T} [\nabla_{v_I} f]^T v_I^0\right)$ [1].

Because both methods investigate physical power flow models, a reduced network with fewer nodes can be constructed. However, if the voltages at *Ibus* (i.e., rotor angle) differ from the voltages used for linearization, the errors in the voltage estimates can be substantial.

*IV.2. Proposed Method Based on the Holomorphic Embedding*

The holomorphic embedding method is an innovative approach to power flow analysis that tries to discover voltage solutions analytically rather than iteratively [23]. In *HE*'s load modeling, constant impedance (*Z*) load models are represented as $I = Y_{bus}v + I_0$ where *I* represents current, $Y_{bus}$ is the nodal admittance matrix, v represent voltages, and $I_0$ is the designated current. The $Y_{bus}$ matrix has been altered and divided into two parts: $Y_{tr}$ and $Y_{sh}$. $Y_{tr}$ has zero row sums meaning that its row sums are all zeros. In contrast, $Y_{sh}$ is a diagonal matrix made up of nodal shunt elements. If $Y_{sh}$ is zero, the germ solution serves no flows or loads in the system. Voltage magnitudes and angles in the germ solution are homogeneous, with a typical magnitude of unity and an angle of zero. It is possible, however, to select a reference voltage magnitude other than unity, which is denoted by $E_{ref}$.

Depending on the bus types specified in [17], the power balance equation is (1A) for Bus $i \in Ibus$ and (1B) for Bus $j \in KM$:

$$e_i^T Y_{tr} v(s) = \frac{1}{E_i^2}[p_i(s) - jq_i(s)]v_i(s) + sI_i - se_i^T Y_{sh} v(s) \tag{1A}$$

$$e_j^T Y_{tr} v(s) = s(p_j - jq_j)w_j^*(s) + sI_j - se_j^T Y_{sh} v(s) \tag{1B}$$

Here, $E_i$ os the voltage magnitude at Bus *i*; $p_i$ and $q_i$ are the real and reactive power injections at Bus *i*, while $p_j$ and $q_j$ are the real and reactive power injections at Bus *j*; and $w_j$ is the reciprocal voltages at Bus *j*, i.e., $w_j v_j = 1$, respectively.

The formulae for $Y_{tr}v$ and $Y_{sh}v$, given the partitioned into submatrices $Y_{tr}$ and $Y_{sh}$ are as follows: $Y_{tr}v = \begin{bmatrix} Y_{tr}^{II} & Y_{tr}^{IKM} \\ Y_{tr}^{KMI} & Y_{tr}^{KMKM} \end{bmatrix} \begin{pmatrix} v_I \\ v_{KM} \end{pmatrix}$ and $Y_{sh}v = \begin{bmatrix} Y_{sh}^I & 0 \\ 0 & Y_{sh}^{KM} \end{bmatrix} \begin{pmatrix} v_I \\ v_{KM} \end{pmatrix}$. According to the germ solution in [23], because the row sums of $Y_{tr}$ are all zero, one can set $v[0]$ or $E_{ref} \angle \theta_0$ to unity. Let $\theta_0 = 0$ or $v[0] = E_{ref}$ for simplicity. According to the definition of w, $v_j(s)w_j(s) = 1$, resulting in: $v_j^*(s)w_j^*(s) = 1$. Consequently, we obtain $w_j^*[0] = \frac{1}{E_{ref}}, w_j^*[1] = -\frac{1}{E_{ref}^2} v_j^*[1]$, and $w_j^*[n] = -\frac{1}{E_{ref}} \sum_{k=0}^{n-1} w_j^*[k] v_j^*[n-k]$ where $j$ and $n$ are indices while $E_{ref}$ is a constant. Assuming that $w_j^*[k]$ and $v_j^*[n-k]$ satisfy the linear equations $w_j^*[k] = A_{w_j^*[k]} \binom{x_I}{y_I} + a_{w_j^*[k]}$ and $v_j^*[n-k] = A_{v_j^*[n-k]} \binom{x_I}{y_I} + a_{v_j^*[n-k]}$ for $k = 0, 1, ..., $ n-1 where $A$ is a row vector and $a$ is a scalar. We can write $w_j^*[n]$ as:

$$w_j^*[n] = -\frac{1}{E_{ref}} \begin{bmatrix} A_{w_j^*[0]} \binom{x_I}{y_I} + a_{w_j^*[0]} & \cdots & A_{w_j^*[n-1]} \binom{x_I}{y_I} + a_{w_j^*[n-1]} \end{bmatrix} \begin{pmatrix} A_{v_j^*[n]} \binom{x_I}{y_I} + a_{v_j^*[n]} \\ \vdots \\ A_{v_j^*[n]} \binom{x_I}{y_I} + a_{v_j^*[1]} \end{pmatrix}$$

$$= \begin{pmatrix} x_I \\ y_I \\ 1 \end{pmatrix}^T M_{w_j^*[n]} \begin{pmatrix} x_I \\ y_I \\ 1 \end{pmatrix} \tag{2A}$$

Where $M_{w_j^*[n]} = \frac{1}{2}(M_{AA[n]} + M_{AA[n]}^T)$ and $M_{AA[n]} = -\frac{1}{E_{ref}} \begin{bmatrix} A_{w_j^*[0]}^T & \cdots & A_{w_j^*[n-1]}^T \\ a_{w_j^*[0]}^T & \cdots & a_{w_j^*[n-1]}^T \end{bmatrix} \begin{bmatrix} A_{v_j^*[n]} & a_{v_j^*[n]} \\ \vdots & \vdots \\ A_{v_j^*[n]} & a_{v_j^*[1]} \end{bmatrix}$.

(2A) is quadratic in *Ibus* voltages, but it can be linearized around known bus voltage values as follows:

$$w_j^*[n] \cong A_{w_j^*[n]} \binom{x_I}{y_I} + a_{w_j^*[n]} \tag{2B}$$

Where $A_{w_j^*[n]} = 2A_{w_j^*[n]}^0$, $a_{w_j^*[n]} = a_{w_j^*[n]}^0 - A_{w_j^*[n]}^0 \binom{x_I^0}{y_I^0}$, and $\begin{bmatrix} A_{w_j^*[n]}^0 & a_{w_j^*[n]}^0 \end{bmatrix} = \begin{pmatrix} x_I^0 \\ y_I^0 \\ 1 \end{pmatrix}^T M_{w_j^*[n]}$.

From the linear expression $v_j^*[n] = A_{v_j^*[n]} \binom{x_I}{y_I} + a_{v_j^*[n]}$ where $A_{v_j[n]} = \left(A_{v_j^*[n]}\right)^*$ and $a_{v_j[n]} = \left(a_{v_j^*[n]}\right)^*$, we find $v_j[n] = A_{v_j[n]} \binom{x_I}{y_I} + a_{v_j[n]}$. With $w_j^*[n]$ linearized as shown in (2B) and $v_I(s) = E_{ref}\mathbf{1} + (v_I - E_{ref}\mathbf{1})s$, (1B) becomes:

$$e_j^T Y_{tr}^{KMKM} v_{KM}[n+1] = e_j^T B_{v_{KM}[n+1]} v_I + e_j^T b_{v_{KM}[n+1]} \tag{2C}$$

Where $e_j^T B_{v_{KM}[n+1]} = -\delta_{n0}e_j^T Y_{tr}^{KMI} + (p_j - \mathbf{j}q_j)A_{w_j^*[n]} - Y_{sh}^j A_{v_j[n]}$ and $e_j^T b_{v_{KM}[n+1]} = \delta_{n0}e_j^T Y_{tr}^{KMI} E_{ref} 1 + (p_j - \mathbf{j}q_j)a_{w_j^*[n]} + I_j\delta_{n1} - Y_{sh}^j a_{v_j[n]}$. (2C) is solved because $Y_{tr}^{KMKM}$ is a full-rank matrix:

$$v_{KM}[n+1] = \tilde{A}_{v_j[n+1]} v_I + \tilde{a}_{v_j[n+1]} \quad \text{or}$$

$$\begin{pmatrix} x_{KM}[n+1] \\ y_{KM}[n+1] \end{pmatrix} = \begin{bmatrix} \tilde{A}_{v_j[n+1]}^{Re} & -\tilde{A}_{v_j[n+1]}^{Im} \\ \tilde{A}_{v_j[n+1]}^{Im} & \tilde{A}_{v_j[n+1]}^{Re} \end{bmatrix} \begin{pmatrix} x_I \\ y_I \end{pmatrix} + \begin{pmatrix} \tilde{a}_{v_j[n+1]}^{Re} \\ \tilde{a}_{v_j[n+1]}^{Im} \end{pmatrix} = A_{v_j[n+1]} \begin{pmatrix} x_I \\ y_I \end{pmatrix} + a_{v_j[n+1]} \quad (2D)$$

Where $A_{v_j[n+1]} = (Y_{tr}^{KMKM})^{-1} B_{v_{KM}[n+1]}$ and $a_{v_j[n+1]} = (Y_{tr}^{KMKM})^{-1} b_{v_{KM}[n+1]}$.

The cardinalities of $A_{v_j[n+1]}$ and $a_{v_j[n+1]}$ are $n_{KM} \times n_I$ and $n_{KM} \times 1$, respectively. Summation over $n$ givess: $\begin{pmatrix} x_{KM} \\ y_{KM} \end{pmatrix} = \left\{ \sum_{n=0}^{\infty} s^n A_{v_j[n]} \right\} \begin{pmatrix} x_I \\ y_I \end{pmatrix} + \left\{ \sum_{n=0}^{\infty} s^n a_{v_j[n]} \right\}$. Because adding an infinite number of terms is infeasible, we seek the Padé approximation [24], which was inspired by [23]:

$$\begin{pmatrix} x_{KM} \\ y_{KM} \end{pmatrix} = \left\{ \sum_{i=0}^{\infty} s^i A_{v_j[i]} \right\} \begin{pmatrix} x_I \\ y_I \end{pmatrix} + \left\{ \sum_{i=0}^{\infty} s^i a_{v_j[i]} \right\} = \frac{\left\{ \sum_{k=0}^{l} s^k C_{v_j[k]} \right\} \begin{pmatrix} x_I \\ y_I \end{pmatrix} + \left\{ \sum_{k=0}^{l} s^k c_{v_j[k]} \right\}}{\sum_{k=0}^{m} s^k b_{v_j[k]}} \quad (3A)$$

$$\rightarrow \left\{ \sum_{k=0}^{l} s^k C_{v_j[k]} \right\} \begin{pmatrix} x_I \\ y_I \end{pmatrix} + \left\{ \sum_{k=0}^{l} s^k c_{v_j[k]} \right\} = \left( \left\{ \sum_{i=0}^{\infty} s^i A_{v_j[i]} \right\} \begin{pmatrix} x_I \\ y_I \end{pmatrix} + \left\{ \sum_{i=0}^{\infty} s^i a_{v_j[i]} \right\} \right) \sum_{k=0}^{m} s^k b_{v_j[k]}$$

In this study, we set $l$ equal to $m$. The first linear system is constructed using the coefficients of $s^k$, $k = 0, \ldots, l$:

For $s^0$: $\left( A_{v_j[0]} \begin{pmatrix} x_I \\ y_I \end{pmatrix} + a_{v_j[0]} \right) b_{v_j[0]} = C_{v_j[0]} \begin{pmatrix} x_I \\ y_I \end{pmatrix} + c_{v_j[0]}$

$$\rightarrow C_{v_j[0]} = b_{v_j[0]} A_{v_j[0]}, c_{v_j[0]} = b_{v_j[0]} a_{v_j[0]} \quad (3B)$$

For $s^1$: $\left( A_{v_j[0]} \begin{pmatrix} x_I \\ y_I \end{pmatrix} + a_{v_j[0]} \right) b_{v_j[1]} + \left( A_{v_j[1]} \begin{pmatrix} x_I \\ y_I \end{pmatrix} + a_{v_j[1]} \right) b_{v_j[0]} = C_{v_j[1]} \begin{pmatrix} x_I \\ y_I \end{pmatrix} + c_{v_j[1]}$

$$\rightarrow C_{v_j[1]} = b_{v_j[1]} A_{v_j[0]} + b_{v_j[0]} A_{v_j[1]}, c_{v_j[1]} = b_{v_j[1]} a_{v_j[0]} + b_{v_j[0]} a_{v_j[1]} \quad (3C)$$

For $s^l$: $\left(A_{v_j[0]}\binom{x_I}{y_I} + a_{v_j[0]}\right)b_{v_j[l]} + \cdots + \left(A_{v_j[l]}\binom{x_I}{y_I} + a_{v_j[l]}\right)b_{v_j[0]} = C_{v_j[l]}\binom{x_I}{y_I} + c_{v_j[l]}$

$$\rightarrow C_{v_j[l]} = \sum_{k=0}^{l} b_{v_j[l-k]}A_{v_j[k]}, c_{v_j[l]} = \sum_{k=0}^{l} b_{v_j[l-k]}a_{v_j[k]} \tag{3D}$$

The second linear system is given by the coefficients of $s^k$, $k = l+1, \ldots, l+m$:

For $s^{l+1}$: $\left(A_{v_j[l+1]}\binom{x_I}{y_I} + a_{v_j[l+1]}\right)b_{v_j[0]} + \cdots + \left(A_{v_j[l-m+1]}\binom{x_I}{y_I} + a_{v_j[l-m+1]}\right)b_{v_j[m]} = 0$

$$\rightarrow \begin{bmatrix} A_{v_j[l+1]}e_1 & \cdots & A_{v_j[l-m+1]}e_1 \\ \vdots & \ddots & \vdots \\ A_{v_j[l+1]}e_{n_I} & \cdots & A_{v_j[l-m+1]}e_{n_I} \\ a_{v_j[l+1]} & \cdots & a_{v_j[l-m+1]} \end{bmatrix}^{[(2NI+1)\times 2NKM]\times(m+1)} \begin{pmatrix} b_{v_j[0]} \\ \vdots \\ b_{v_j[m]} \end{pmatrix} = 0 \tag{3E}$$

For $s^{l+2}$: $\left(A_{v_j[l+2]}\binom{x_I}{y_I} + a_{v_j[l+2]}\right)b_{v_j[0]} + \cdots + \left(A_{v_j[l-m+2]}\binom{x_I}{y_I} + a_{v_j[l-m+2]}\right)b_{v_j[m]} = 0$

$$\rightarrow \begin{bmatrix} A_{v_j[l+2]}e_1 & \cdots & A_{v_j[l-m+2]}e_1 \\ \vdots & \ddots & \vdots \\ A_{v_j[l+2]}e_{n_I} & \cdots & A_{v_j[l-m+2]}e_{n_I} \\ a_{v_j[l+2]} & \cdots & a_{v_j[l-m+2]} \end{bmatrix}^{[(2NI+1)\times 2NKM]\times(m+1)} \begin{pmatrix} b_{v_j[0]} \\ \vdots \\ b_{v_j[m]} \end{pmatrix} = 0 \tag{3F}$$

For $s^{l+m}$: $\left(A_{v_j[l+m]}\binom{x_I}{y_I} + a_{v_j[l+m]}\right)b_{v_j[0]} + \cdots + \left(A_{v_j[l]}\binom{x_I}{y_I} + a_{v_j[l]}\right)b_{v_j[m]} = 0$

$$\rightarrow \begin{bmatrix} A_{v_j[l+m]}e_1 & \cdots & A_{v_j[l]}e_1 \\ \vdots & \ddots & \vdots \\ A_{v_j[l+m]}e_{n_I} & \cdots & A_{v_j[l]}e_{n_I} \\ a_{v_j[l+m]} & \cdots & a_{v_j[l]} \end{bmatrix}^{[(2NI+1)\times 2NKM]\times(m+1)} \begin{pmatrix} b_{v_j[0]} \\ \vdots \\ b_{v_j[m]} \end{pmatrix} = 0 \tag{3G}$$

It is worth noticing that the linear matrix equations in (3E)-(3G) are derived in such a way that the superposition $s$ becomes zero at any value for $x_I$ and $y_I$. To find a linear equation, use (3E)-(3G):

$$W_b^{[m(2NI+1)\times 2NKM]\times(m+1)} \begin{pmatrix} b_{v_j[0]} \\ \vdots \\ b_{v_j[m]} \end{pmatrix} = 0 \text{ where } W_b = \begin{bmatrix} A_{v_j[l+1]}e_1 & \cdots & A_{v_j[l-m+1]}e_1 \\ \vdots & \ddots & \vdots \\ A_{v_j[l+m]}e_{n_I} & \cdots & A_{v_j[l]}e_{n_I} \\ a_{v_j[l+1]} & \cdots & a_{v_j[l-m+1]} \\ \vdots & \ddots & \vdots \\ a_{v_j[l+m]} & \cdots & a_{v_j[l]} \end{bmatrix} \tag{4A}$$

As a result, the vector formed of the coefficients in the denominator series of the Padé approximation resides in the null space of the $W_b$ matrix. Because the denominator function of the Padé approximation (3A) cannot be zero, the vector must be a nontrivial solution. To explore the matrix's null space, the singular value decomposition of $W_b$ ($W_b = U\Sigma V^T$) is used:

$$\begin{pmatrix} b_{v_j[0]} \\ \vdots \\ b_{v_j[m]} \end{pmatrix} = V^T e_M \qquad (4B)$$

The coefficients in the nominator series are found utilizing the values for $b_{v_j[k]}$ from (3B)-(3D), and (4B). The Padé approximation yields a linear relationship between $v_I$ and $v_{KM}$ at $s = 1$:

$$\begin{pmatrix} x_{KM} \\ y_{KM} \end{pmatrix} = \left\{ \sum_{i=0}^{\infty} A_{v_j[i]} \right\} \begin{pmatrix} x_I \\ y_I \end{pmatrix} + \left\{ \sum_{i=0}^{\infty} a_{v_j[i]} \right\} = H_{KM} \begin{pmatrix} x_I \\ y_I \end{pmatrix} + h_{KM} \qquad (5)$$

Where $H_{KM} = \frac{\sum_{k=0}^{l} C_{v_j[k]}}{\sum_{k=0}^{m} b_{v_j[k]}}$ and $h_{KM} = \frac{\sum_{k=0}^{l} c_{v_j[k]}}{\sum_{k=0}^{m} b_{v_j[k]}}$.

The power balance equation for *Ibus* is $p_i(s) - \mathbf{j}q_i(s) = v_i^*(s)[e_i^T Y_{tr} v(s) - sI_i + sy_{sh}^i v_i(s)]$. Linear equations that $v_i^*(s)[e_i^T Y_{tr} v(s) - sI_i + sy_{sh}^i v_i(s)] = \tilde{F}_i(s)\begin{pmatrix} x_I \\ y_I \end{pmatrix} + \tilde{f}_i(s)$ are obtained by linearizing the right-hand side of the power balance equation at *Ibus* using (5) and $v_I(s) = E_{ref} + (v_I - E_{ref})s$. As a result,

$$p_i(s) - \mathbf{j}q_i(s) = \tilde{F}_i(s)\begin{pmatrix} x_I \\ y_I \end{pmatrix} + \tilde{f}_i(s) \rightarrow \begin{cases} p_I[k] = \tilde{F}^{Re}[k]\begin{pmatrix} x_I \\ y_I \end{pmatrix} + \tilde{f}^{Re}[k] \\ q_I[k] = -\tilde{F}^{Im}[k]\begin{pmatrix} x_I \\ y_I \end{pmatrix} - \tilde{f}^{Im}[k] \end{cases} \qquad (7A)$$

As with (3A)-(3I), the Padé approximation is applied to produce the linear expression for real and reactive power generation at *Ibus*:

$$p_I = H_P \begin{pmatrix} x_I \\ y_I \end{pmatrix} + h_P \text{ and } q_I = H_Q \begin{pmatrix} x_I \\ y_I \end{pmatrix} + h_Q \qquad (7B)$$

The definitions of $\xi_i \ (\equiv q_i x_i)$ and $\zeta_i \ (\equiv q_i y_i)$ reveal the following:

$$\xi_i[n] = \sum_{k=0}^{n} q_i[n-k]x_i[k]$$

$$= E_{ref}\left\{e_i^T H_{Q[n]}\begin{pmatrix}x_I\\y_I\end{pmatrix} + e_i^T h_{Q[n]}\right\} + (x_i - E_{ref})\left\{e_i^T H_{Q[n-1]}\begin{pmatrix}x_I\\y_I\end{pmatrix} + e_i^T h_{Q[n-1]}\right\}$$

$$= e_i^T E_{ref} H_{Q[n]}\begin{pmatrix}x_I\\y_I\end{pmatrix} - e_i^T E_{ref} H_{Q[n-1]}\begin{pmatrix}x_I\\y_I\end{pmatrix} + e_i^T h_{Q[n-1]}[e_i^T \quad 0]\begin{pmatrix}x_I\\y_I\end{pmatrix}$$

$$+ \begin{pmatrix}x_I\\y_I\end{pmatrix}^T H_{Q[n-1]}^I \begin{pmatrix}x_I\\y_I\end{pmatrix} + e_i^T E_{ref} h_{Q[n]} - e_i^T E_{ref} h_{Q[n-1]}$$

$$\cong \left\{e_i^T E_{ref} H_{Q[n]} - e_i^T E_{ref} H_{Q[n-1]} + e_i^T h_{Q[n-1]}[e_i^T \quad 0] + 2\begin{pmatrix}x_I^0\\y_I^0\end{pmatrix}^T H_{Q[n-1]}^I\right\}\begin{pmatrix}x_I\\y_I\end{pmatrix}$$

$$+ e_i^T E_{ref} h_{Q[n]} - e_i^T E_{ref} h_{Q[n-1]} - \begin{pmatrix}x_I^0\\y_I^0\end{pmatrix}^T H_{Q[n-1]}^I \begin{pmatrix}x_I^0\\y_I^0\end{pmatrix}$$

$$= e_i^T H_{\xi[n]}\begin{pmatrix}x_I\\y_I\end{pmatrix} + e_i^T h_{\xi[n]} \tag{8A}$$

$$\zeta_i[n] = \sum_{k=0}^{n} q_i[k]y_i[n-k] = y_i\left\{e_i^T H_{Q[n-1]}\begin{pmatrix}x_I\\y_I\end{pmatrix} + e_i^T h_{Q[n-1]}\right\}$$

$$= \begin{pmatrix}x_I\\y_I\end{pmatrix}^T H_{Q[n-1]}^{II}\begin{pmatrix}x_I\\y_I\end{pmatrix} + e_i^T h_{Q[n-1]}[0 \quad e_i^T]\begin{pmatrix}x_I\\y_I\end{pmatrix}$$

$$\cong \left\{2\begin{pmatrix}x_I^0\\y_I^0\end{pmatrix}^T H_{Q[n-1]}^{II} + e_i^T h_{Q[n-1]}[0 \quad e_i^T]\right\}\begin{pmatrix}x_I\\y_I\end{pmatrix} - \begin{pmatrix}x_I^0\\y_I^0\end{pmatrix}^T H_{Q[n-1]}^{II}\begin{pmatrix}x_I^0\\y_I^0\end{pmatrix}$$

$$= e_i^T H_{\zeta[n]}\begin{pmatrix}x_I\\y_I\end{pmatrix} + e_i^T h_{\zeta[n]} \tag{8B}$$

Where $H_{Q[k]}^I = \frac{1}{2}\left\{\begin{bmatrix}e_j\\0\end{bmatrix}e_i^T H_{Q[k]} + H_{Q[k]}^T e_j[e_i^T \quad 0]\right\}$ and $H_{Q[k]}^{II} = \frac{1}{2}\left\{\begin{bmatrix}0\\e_j\end{bmatrix}e_i^T H_{Q[k]} + H_{Q[k]}^T e_j[0 \quad e_i^T]\right\}$.

Like (3A)-(3I), the Padé approximation offers the linear formula for $\xi_I$ and $\zeta_I$ at *Ibus*:

$$\xi_I = H_\xi\begin{pmatrix}x_I\\y_I\end{pmatrix} + h_\xi \text{ and } \zeta_I = H_\zeta\begin{pmatrix}x_I\\y_I\end{pmatrix} + h_\zeta \tag{9}$$

The linearization results are summarized below:

$$\chi \equiv \begin{pmatrix} x_{KM} \\ y_{KM} \\ p_I \\ q_I \\ \xi_I \\ \zeta_I \end{pmatrix} = \begin{bmatrix} H_{KM} \\ H_P \\ H_Q \\ H_\xi \\ H_\zeta \end{bmatrix} \begin{pmatrix} x_I \\ y_I \end{pmatrix} + \begin{pmatrix} h_{KM} \\ h_P \\ h_Q \\ h_\xi \\ h_\zeta \end{pmatrix} = H \begin{pmatrix} x_I \\ y_I \end{pmatrix} + h \quad (10)$$

*IV.3. Error Analysis*

The number of terms in the nominator and denominator series, denoted by $l$ and $m$, respectively, influences the accuracy of the Padé approximation. While not as robust as the McLaurin series, the amount of terms utilized has an effect on the precision of the Padé approximation [30]. The holomorphic embedding procedure limits the series to a limited number of terms, resulting in errors in the matrices and vectors due to inaccuracies in the data and observations. The errors are caused by two factors: deviation from the voltage values at *Ibus* during the linearization process and truncation error due to the finite nature of *l*. The matrix and currents obtained using the Padé approximation are represented by $\overline{\Psi}_l$ and $\overline{\psi}_l$, respectively. In contrast, $\Psi_\infty$ and $\psi_\infty$ denote the matrices and currents acquired by employing an infinite number of terms. Denote $\begin{pmatrix} \chi_{true} \\ v_I^{true} \end{pmatrix}$ and $\begin{pmatrix} \chi_{LS} \\ v_I^{LS} \end{pmatrix}$ as the solutions with and without errors in (10): $\begin{pmatrix} \chi_c^{LS} \\ v_I^{LS} \end{pmatrix} = \underset{\chi, v_I}{\mathrm{argmin}} \left\| \overline{\psi}_l - \overline{\Psi}_l \begin{pmatrix} \chi \\ v_I \end{pmatrix} \right\|$ and $\begin{pmatrix} \chi_{true} \\ v_I^{true} \end{pmatrix} = \underset{\chi, v_I}{\mathrm{argmin}} \left\| \psi_\infty - \Psi_\infty \begin{pmatrix} \chi \\ v_I \end{pmatrix} \right\|$. The difference between solutions is bounded by [31]:

$$\frac{\left\| \begin{pmatrix} \chi_c^{LS} \\ v_I^{LS} \end{pmatrix} - \begin{pmatrix} \chi_{true} \\ v_I^{true} \end{pmatrix} \right\|_2}{\left\| \begin{pmatrix} \chi_c^{LS} \\ v_I^{LS} \end{pmatrix} \right\|_2} \leq \epsilon \{ 2 \cos^{-1}\theta \kappa_2(\overline{\Psi}_l) + \tan\theta \kappa_2(\overline{\Psi}_l)^2 \} \quad (11)$$

Where $\quad \epsilon = max\left( \frac{\|\Psi_\infty - \overline{\Psi}_l\|_2}{\|\overline{\Psi}_l\|_2}, \frac{\|\psi_\infty - \overline{\psi}_L\|_2}{\|\overline{\psi}_L\|_2} \right) \quad , \quad \theta = \sin^{-1}\left( \frac{\left\| \overline{\psi}_l - \overline{\Psi}_l \begin{pmatrix} \chi_c^{LS} \\ v_I^{LS} \end{pmatrix} \right\|_2}{\|\overline{\psi}_l\|_2} \right) \quad ,$

$\kappa_2(\overline{\Psi}_l) = \|\overline{\Psi}_l\|_2 \left\| (\overline{\Psi}_l^T \overline{\Psi}_l)^{-1} \overline{\Psi}_l \right\|_2$, and $\kappa_2(\overline{\Psi}_l)^2 = \|\overline{\Psi}_l\|_2^2 \left\| (\overline{\Psi}_l^T \overline{\Psi}_l)^{-1} \right\|_2$.

(11) is an equality constrained least squares problem with a fixed voltage reference angle, as given in [31]:

$\begin{pmatrix} \chi_c^{LS} \\ v_I^{LS} \end{pmatrix} = \underset{\chi, v_I}{\operatorname{argmin}} \left\| \bar{\psi}_l - \bar{\Psi}_l \begin{pmatrix} \chi \\ v_I \end{pmatrix} \right\| : s.t. \, e_{ref}^T \begin{pmatrix} \chi \\ v_I \end{pmatrix} = 0$. The elementary column vector corresponding to the location of the reference angle node (ref) is denoted by $e_{ref}$. Let $\widetilde{\Psi}_l$ be the partitioned matrix $\bar{\Psi}_l$ from which the node $ref$ is removed. The solution to the equality constrained least squares problem is [31]: $\begin{pmatrix} \chi_c^{LS} \\ v_I^{LS} \end{pmatrix}_{ref} = \left( \widetilde{\Psi}_l^T \widetilde{\Psi}_l \right)^{-1} \widetilde{\Psi}_l^T \bar{\psi}_l$. The subscript $ref$ indicates that the element in $\chi_c^{LS}$ corresponding to the reference angle node has been removed. To complete $\begin{pmatrix} \chi_c^{LS} \\ v_I^{LS} \end{pmatrix}$, add one row to the element in $\begin{pmatrix} \chi_c^{LS} \\ v_I^{LS} \end{pmatrix}_{ref}$ by including a row vector and an element on the right-hand side in $\left( \widetilde{\Psi}_l^T \widetilde{\Psi}_l \right)^{-1} \widetilde{\Psi}_l^T \bar{\psi}_l$. As previously stated, (10) is only applicable in the neighborhood of the dispatch where linearization is performed. The linear coefficients in (10) must be revised if the voltages at *Ibus* change significantly.

## V. ANALYTICAL SOLUTION TO SWING EQUATIONS

When there is a disturbance in a power system, the swing equation governs the system dynamics: $M_j \frac{d^2 \delta_j}{dt^2} + D_j \frac{d\delta_j}{dt} = p_j^{mech} - p_j^{elec}$ where $M_j$ and $D_j$ are the inertia and damping coefficients of the $j^{th}$ generator, respectively; $\delta_j$, $p_j^{mech}$, and $p_j^{elec}$ are the rotor angle, mechanical power input, and electrical power output of the $j^{th}$ generator, respectively. Kirchhoff's laws govern the output of electrical power. The coordinate transformation of the swing equation to [17] yields:

$$M_j \frac{d^2 x_j}{dt^2} + D_j \frac{dx_j}{dt} + b_{jj}^0 E_j^2 x_j + \left( p_j^{mech} - g_{jj} E_j^2 \right) y_j - \xi_j - |\tilde{y}_{jk}| E_j^2 v_k^x = 0$$

$$M_j \frac{d^2 y_j}{dt^2} + D_j \frac{dy_j}{dt} + b_{jj}^0 E_j^2 y_j - (p_j^{mech} - g_{jj} E_j^2) x_j - \zeta_j - |\tilde{y}_{jk}| E_j^2 v_k^y = 0 \tag{12}$$

where $b_{jj}^0 = b_{jj} + \frac{M_j}{E_j^2}(d\delta_j/dt)^2$.

Using (10), (12) becomes:

$$M_j e_j^T \frac{d^2 w_I}{dt^2} + D_j e_j^T \frac{dw_I}{dt} + \{b_{jj}^0 E_j^2 e_j^T + (p_j^{mech} - g_{jj} E_j^2) e_{NI+j}^T - e_j^T H_\xi - |\tilde{y}_{jk}| E_j^2 e_k^T H_{KM}\} w_I$$

$$= e_j^T h_\xi + |\tilde{y}_{jk}| E_j^2 e_k^T h_{KM} \tag{13}$$

$$M_j e_{NI+j}^T \frac{d^2 w_I}{dt^2} + D_j e_{NI+j}^T \frac{dw_I}{dt} + \{b_{jj}^0 E_j^2 e_{NI+j}^T - (p_j^{mech} - g_{jj} E_j^2) e_j^T - e_j^T H_\zeta - |\tilde{y}_{jk}| E_j^2 e_{NKM+k}^T H_{KM}\} w_I$$

$$= e_j^T h_\zeta + |\tilde{y}_{jk}| E_j^2 e_{NKM+k}^T h_{KM} \tag{14}$$

Where $w_I = (x_I^T \quad y_I^T)^T$. It should be noted that the factors in square brackets in (13) and (14) are related to the admittance matrix and are normalized to the inertia of the corresponding generators. (13) and (14) are both simplified [17]:

$$\frac{d^2 w_I}{dt^2} + diag\left(\frac{D_I}{M_I}\right) \frac{dw_I}{dt} + L w_I + l = 0 \tag{15}$$

Where $e_j^T L = \frac{1}{M_j} [b_{jj}^0 e_j^T + (p_j^{mech} - g_{jj} E_j^2) e_{NI+j}^T - e_j^T H_\xi - |\tilde{y}_{jk}| E_j^2 e_k^T H_{KM}]$, $e_{NI+j}^T L = \frac{1}{M_j} [b_{jj}^0 e_{NI+j}^T - (p_j^{mech} - g_{jj} E_j^2) e_j^T - e_j^T H_\zeta - |\tilde{y}_{jk}| E_j^2 e_{NK+k}^T H_K]$, $e_j^T l = -\frac{1}{M_j}(e_j^T h_\xi + |\tilde{y}_{jk}| E_j^2 e_k^T h_K)$, and $e_{NI+j}^T l = -\frac{1}{M_j}(e_j^T h_\zeta + |\tilde{y}_{jk}| E_j^2 e_{NK+k}^T h_K)$, respectively.

The linear differential equation (15), like the one in [17], is constrained by two parameters, O1 and O2. Because of these constraints, (15) is equivalent to the swing equation (11). O1 denotes all of the voltages at *Ibus*, whereas O2 denotes the boundary condition where the linear approximation (10) is still valid. (17) depicts the combination of O1 and O2:

$$O1: \left\|\left(1 - \frac{1}{E_j}\sqrt{x_j^2 + y_j^2}; x_j \frac{dx_j}{dt} + y_j \frac{dy_j}{dt}\right)\right\| = 0 \text{ and } O2: \left|\frac{\|w_I\|}{\|w_I^0\|} - 1\right| \leq \varepsilon_{O2} \quad (16)$$

$$\left\|\left(1 - \frac{1}{E_j}\sqrt{x_j^2 + y_j^2}; x_j \frac{dx_j}{dt} + y_j \frac{dy_j}{dt}; \frac{p_G}{p_G^0} - 1\right)\right\| \leq \varepsilon \quad (17)$$

The validity region is defined by (17), as indicated in [17]. When (17) is satisfied, (15) remains within a distance of $\varepsilon$ from the solution to the swing equation. As a result, (15) is also $\varepsilon$-close to the true solution of the swing equation, and the error is bounded by $k_{(17)}\varepsilon$, where $k_{(17)}$ is a constant as proved in [17]. According to [32], velocity control in the form of $x_j \frac{dx_j}{dt} + y_j \frac{dy_j}{dt} = 0$ is more effective in maintaining O1 than location control. The following is the solution to differential equation (15):

$$\begin{pmatrix} w_I \\ dw_I/dt \end{pmatrix} = \sum_{k=1}^{4NI} \beta_k^* \Psi_k u + \begin{pmatrix} -L^{-1}l \\ 0 \end{pmatrix} \quad (18)$$

Where, as demonstrated in [17], $\Psi_k$ is the matrix form of the $k^{th}$ vector spanning the null space of $I \otimes T - D_l^T \otimes T$ where $T = \begin{bmatrix} 0 & I \\ -L & -diag\left(\frac{D_I}{M_I}\right) \end{bmatrix}$ and $D_l$ is a block-diagonal matrix containing the eigenvalues of $T$; $u = \begin{pmatrix} u_{Re} \\ u_{Co} \end{pmatrix}$ where $u_{Re} = \begin{bmatrix} exp(\lambda_{Re}^1 t) \\ \vdots \\ exp(\lambda_{Re}^{N_{Re}} t) \end{bmatrix}$ and $u_{Re} = \begin{bmatrix} exp(\lambda_{Co}^1 t)cos(\lambda_{Co}^2 t) \\ exp(\lambda_{Co}^1 t)sin(\lambda_{Co}^2 t) \\ \vdots \\ exp(\lambda_{Co}^{N_{Co}-1} t)cos(\lambda_{Co}^{N_{Co}} t) \\ exp(\lambda_{Co}^{N_{Co}-1} t)sin(\lambda_{Co}^{N_{Co}} t) \end{bmatrix}$ ; and $\beta^* = \underset{\beta}{\text{argmin}} \left\|[\psi_1 \cdots \psi_{4NI}]\beta - \begin{pmatrix} v_I^0 + L^{-1}l \\ -Jv_I^0 \end{pmatrix}\right\|$ where $\psi_l = \Psi_l u(t=0)$ and $J = \begin{bmatrix} 0 & \omega_j^0 \\ -\omega_j^0 & 0 \end{bmatrix}$.

The rotor angle $\delta_I$ is determined as the imaginary component of $log(x_I + jy_I)$. Because the reference angle bus is fixed while calculating the sensitivities in (10), the rotor angles are relative to the reference angle node ref. To put it another way, the estimated rotor angles are $\delta_I^{(18)} = \delta_I - \delta_{ref}$. $\delta_I$ is frequently referenced against the center of inertia $\delta_{COI}$, $\delta_{COI} = \sum_j \left(\frac{M_j}{M_{tot}}\right) \delta_j$ where $M_{tot} = \sum_i M_i$. Using the definition of $\delta_{COI}$, it is possible to determine:

$$\delta_{COI} - \delta_{ref} = \sum_j \left(\frac{M_j}{M_{tot}}\right)(\delta_j - \delta_{ref}) \tag{19}$$

(19) yields (20), which alters the reference angle:

$$\delta_j - \delta_{COI} = \delta_j^{(18)} - \sum_i \left(\frac{M_i}{M_{tot}}\right)\delta_i^{(18)} = \sum_i \left[\delta_{ij} - \frac{M_i}{M_{tot}}\right]\delta_i^{(18)} \tag{20}$$

$\delta_{ij}$ is the Kronecker delta [33]. The rotor angles relative to the moment of inertia can be determined directly using the solution in (18) since the linear coefficients, $\delta_{ij} - \frac{M_i}{M_{tot}}$, are independent of the reference bus. As a result, the linear approximation using $HE$ in (10) remains valid, despite the fact that it is predicted on a specific choice of the reference angle bus.

## VI. BEYOND THE VALIDITY REGION

Denote $f$ and $g$ as (15) and (16), respectively:

$$f(w_I, \omega_I, \rho_I) = \rho_I + diag\left(\frac{D_I}{M_I}\right)\omega_I + Lw_I + l = 0 \tag{21}$$

$$g(v_I, \omega_I) = \begin{bmatrix} E_I^2 - vec(e_j e_j^T + e_{j+N} e_{j+N}^T)^T \omega_I \otimes \omega_I \\ vec(e_j e_j^T + e_{j+N} e_{j+N}^T)^T w_I \otimes \omega_I \end{bmatrix} = 0 \tag{22}$$

$F(w_I, \omega_I, \rho_I)$ is formed as follows:

$$F(v_I, \omega_I, \rho_I) = \begin{bmatrix} f(v_I, \omega_I, \rho_I) \\ g(v_I, \omega_I) \end{bmatrix} = 0 \tag{23}$$

According to the definition of the validity region in [17], the value of $F$ is non-zero at the last point on the boundary. Iteratively, the Newton-Raphson method is used to obtain a consistent initial point for $w_I, \omega_I \left(= \frac{dw_I}{dt}\right)$, and $\rho_I \left(= \frac{d\omega_I}{dt}\right)$. The starting point needs to be compatible with the differential equation (15). Kirchhoff's laws are utilized to compute the initial power injection.

The solution in (18) is valid if the power injection is close enough to the estimated value from (10) and (18) as well as the update for *l*. Otherwise, the linear sensitivities in (10) must be recalculated to determine matrices *H* and *L*, and the solution in (18) must be recalculated, resulting in greater complexity.

## VII. NUMERICAL SIMULATION METHODS

*VII.1. Exact Approach: Time-Domain Simulation (TDS)*

The power flow equation becomes linear when dealing with linear loads such as constant impedance load, allowing for relatively simple digital simulation. The power flow equation becomes quadratic when a constant current or constant power load is present. Several iterative steps are involved in the modified Euler method:

1) Predictor stage: The voltages at *KM* are estimated in this stage by solving the *QPF* problem using the voltages at *Ibus*, $v_I$. Using the given equations $M_j \frac{d\omega_j}{dt} = p_j^{mech} - p_j^{elec} - D_j(\delta_j - \delta_j^0)$ and $\frac{d\delta_j}{dt} = \omega_j - \omega_j^0$, the mechanical power, $p_j^{mech}$, electrical power, $p_j^{elec}$, and damping factor, $D_j$, of each node are utilized to compute the first derivative of angular frequency, $d\omega_j/dt$, and the first derivative of the voltage angle, $d\delta_j/dt$. Using these first-order derivatives, the second-order derivatives of rotor angle, $d\omega_i^{2)}/dt$ and the accompanying frequency and angle values, $d\delta_i^{2)}/dt$, $\omega_i^{2)}$ and $\delta_i^{2)}$, are then calculated.

2) Corrector stage: A new *QPF* problem is solved using the rotor angle values calculated in the predictor stage, $\delta_i^{2)}$. This entails calculating the power generated at each node, which is then used to determine the second-order derivatives of rotor angle, $d\omega_i^{3)}/dt$, and their corresponding frequency and rotor angle values, $d\delta_i^{3)}/dt$, $\omega_i^{3)}$ and $\delta_i^{3)}$, respectively. Solving the *QPF* problem determines power generation.

3) Return to Step 1 for the next iteration based on the revised voltage angles at *Ibus* $\delta_i^{3)}$.

In terms of node categorization, the *QPF* problem at Step 1 differs from conventional *PF* problems. There are reference nodes, *PV* nodes, and *PQ* nodes in conventional *PF* problems. In the *QPF* problem, the voltage magnitudes and angles are known at all *Ibus* while the real and reactive power injections are known at the remaining nodes, which are referred to as *KM*. Normally, the voltage angle at the reference angle bus is set to zero.

The first step of the *QPF* problem has *2N* equations and *2N* unknowns that can be utilized to compute the vector $[v_x; v_y; p_{inj}; q_{inj}]^{4N}$, where $v_x$ and $v_y$ are the real and the imaginary components of the voltage vector, and $p_{inj}$ and $q_{inj}$ are the real and the reactive power injections, respectively. Despite the fact that Step 1 differs from the conventional *PF* problem in terms of node classification, it can still be solved using a method such as the Kronecker product [22]. The Kronecker product is a matrix operation that can convert power flow equations into matrices, making numerical computation easier to solve the equations. The *2N* equations and *2N* unknowns in Step 1's *QPF* problem can be written in matrix form using the Kronecker product and solved using a variety of numerical approaches.

$$\begin{cases} e_j^T v = v_j, & j \in Ibus \\ v^T M_i v = m_i, & i \in K \cup M \end{cases} \leftrightarrow \left\{ vec \begin{bmatrix} 0 & e_j \\ e_j^T & 0 \end{bmatrix} \; vec \begin{bmatrix} M_i & 0 \\ 0 & 0 \end{bmatrix} \right\}^T \times \left[ \begin{pmatrix} v \\ 1 \end{pmatrix} \otimes \begin{pmatrix} v \\ 1 \end{pmatrix} \right] = \begin{pmatrix} 2v_j \\ m_i \end{pmatrix} \quad (24)$$

However, the computational costs of numerically solving (24) increase considerably with each iteration. For example, $10^4$ instances of (24) must be solved to complete a 10-second event with a time step of $10^{-3}$ seconds.

*VII.2. Heuristic Truncated McLaurin Series (HTMS) Approach*

Section *VI.1* emphasizes that the *QPF* problem is structured in the same way as conventional *PF* problems. To investigate this similarity, the *HE* [25] is employed for numerical simulations. One may begin by assuming that the time-dependent rotor angles for both $\delta_j$ and $\omega_j$ as a McLaurin series near $t = 0$: $\delta_j = \sum_{k=0}^{\infty} \delta_j[k]t^k$ as well as $\omega_j = \sum_{k=0}^{\infty} \omega_j[k]t^k$. The series is then truncated to the first $m$ terms.

$$\delta_j = \sum_{k=0}^{m-1} \delta_j[k]t^k \text{ and } \omega_j = \sum_{k=0}^{m-1} \omega_j[k]t^k \text{ for } j \in Ibus \tag{25}$$

With (25), one can deduce expressions for $d\omega_j/dt = \sum_{k=0}^{m-1}(k+1)\omega_j[k+1]t^k$ and $d\delta_j/dt = \sum_{k=0}^{m-1}(k+1)\delta_j[k+1]t^k$. Comparing the $\tau^{th}$ power of $t$ in (25) and the definitions of $\delta$ and $\omega$ leads to (26) that establishes the relationship between two consecutive terms, $\tau$ and $\tau + 1$:

$$\begin{bmatrix} M_j & D_j \\ 0 & 1 \end{bmatrix} \begin{pmatrix} \omega_j[\tau+1] \\ \delta_j[\tau+1] \end{pmatrix} = \frac{1}{\tau+1} \begin{pmatrix} p_j^{mech}[\tau] - p_j^{elec}[\tau] \\ \omega_j[\tau] \end{pmatrix}$$

$$\leftrightarrow \begin{pmatrix} \omega_j[\tau+1] \\ \delta_j[\tau+1] \end{pmatrix} = \frac{1}{(\tau+1)M_j} \begin{bmatrix} 1 & -D_j \\ 0 & M_j \end{bmatrix} \begin{pmatrix} p_j^{mech}[\tau] - p_j^{elec}[\tau] \\ \omega_j[\tau] \end{pmatrix} \text{ for } \tau = 0, \cdots, m-1 \tag{26}$$

At each term, the rotor angles are computed using (26), and two nodal variables, either nodal voltages or power injections, are determined based on the bus type, resulting in a *QPF* problem.

The numerical simulation process begins with a simulation of the system at $\tau = 0$, $\omega_j[0]$ and $\delta_j[0]$ are calculated. The *HE* algorithm uses these values to find an analytic solution to identify voltages and power injections $p_j^{elec}[0]$ at $\tau = 0$. The values for $p_j^{elec}[0]$ and (26) allow us to calculate $\omega_j[1]$ and $\delta_j[1]$ for the next *HE* at $\tau = 1$, and so on.

In (25), the state and algebraic variables $u(t)$ are assumed in a truncated McLaurin series form; $u(t) = \sum_{k=0}^{m} u[k]t^k$. The series' convergence radius is:

$$\left| \frac{u[k+1]t^{k+1}}{u[k]t^k} \right| < 1 \rightarrow t < \left| \frac{u[k]}{u[k+1]} \right| \text{ for all } k \geq k_{min} \tag{27}$$

Because the series used in *HE* are numerically evaluated, there is no theoretical limit to the value of $k_{min}$ necessary to satisfy (28A) as long as $k$ is greater than $k_{min}$. As a result, the number of truncations is heuristically chosen heuristically, despite the fact that the number is important to the accuracy of the approximation. For example, the state and algebraic variables $x$ and $y$ are written as $\sum_{j=1}^{4NI} c_j exp\left[\left(\lambda_j^{Re} + \mathbf{j}\lambda_j^{Re}\right)t\right]$ [17]. The resulting rotor angle (voltage angles at *Ibus*) is given by $tan^{-1}\left\{\frac{\sum_{j=1}^{4NI} c_j exp\left[\left(\lambda_j^{Re}+\mathbf{j}\lambda_j^{Re}\right)t\right]}{\sum_{j=1}^{4NI} d_j exp\left[\left(\lambda_j^{Re}+\mathbf{j}\lambda_j^{Re}\right)t\right]}\right\}$. The time dependence is clearly not polynomial.

Considering two extreme cases, rapidly diverging and rapidly converging, we observe that the function is more accurately approximated as the number of terms *m* in the McLaurin series increases. See Figure 2 that depicts the truncated McLaurin series of *exp(t)cos(3t)* (exponentially diverging case) and *exp(-t)cos(3t)* (exponentially converging case). The convergence region for the truncated McLaurin series for $\delta_j, \omega_j, x_i,$ , and $y_i$ should be the minimum of all $k_{min}$ values defined in (27) at the same time, as given by (28) to ensure proper convergence for all variables:

$$t < min\left\{\left|\frac{\delta_j[m-1]}{\delta_j[m]}\right|, \left|\frac{\omega_j[m-1]}{\omega_j[m]}\right|, \left|\frac{x_i[m-1]}{x_i[m]}\right|, \left|\frac{y_i[m-1]}{y_i[m]}\right|\right\} \text{ for all } i, j, \text{ and } m \qquad (28)$$

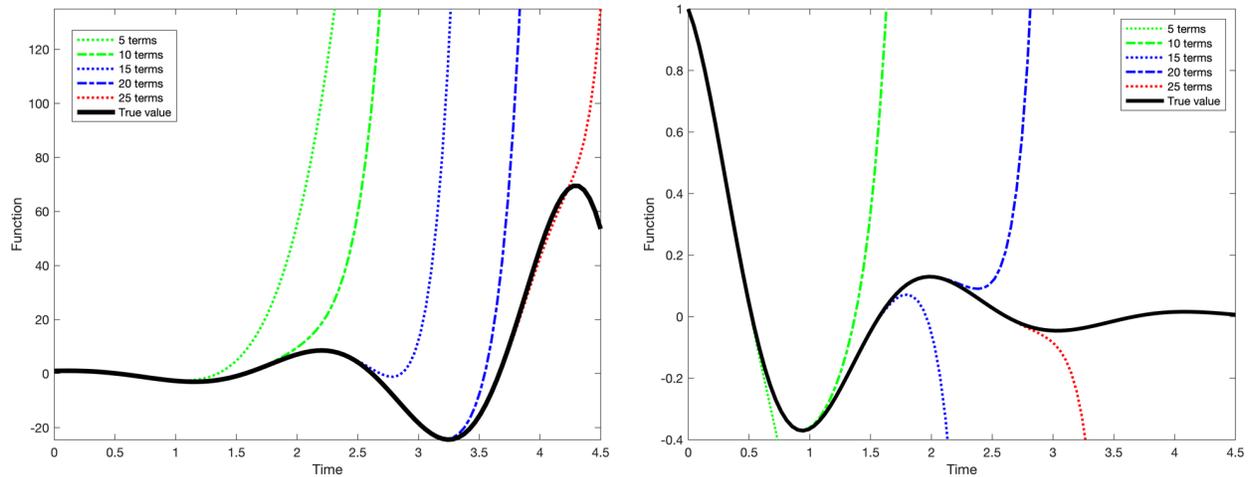

Figure 2. Truncated McLaurin series expansions of an exponentially divulging case (left) and of an exponentially converging case (right). The color lines correspond to the number of terms used for the approximation.

Because the *HE* method is based on numerical term-by-term computation, deriving the analytical forms of $\left|\frac{\delta_j[m-1]}{\delta_j[m]}\right|, \left|\frac{\omega_j[m-1]}{\omega_j[m]}\right|, \left|\frac{x_i[m-1]}{x_i[m]}\right|$, and $\left|\frac{y_i[m-1]}{y_i[m]}\right|$ is difficult, if attainable. To the best of our knowledge, predictinging *t* (corresponding to *m*) in advance to satisfy (28) is thus impossible. When the truncation error is too significant, it can be decreased by either 1) increasing the number of terms (i.e., increasing *m*) or 2) using the Taylor series expansion when the truncation error becomes too significant, such as $\delta_j = \sum_{k=0}^{m-1} \delta_j[k](t - t_{tr})^k$.

To keep the error within an acceptable range, *m* may need to increase exponentially as the simulation's termination time $t_{te}$ increases, As a result, the second approach is employed in this study. Because this method does not involve updating the state and algebraic variables at each iteration, its computation cost is lower than that of traditional time-domain simulation methods. The accuracy of this approach, however, may be jeopardized if the system evolves to a new equilibrium in which the McLaurin series no longer accurately approximates the system dynamics. In such instances as shown in [25], , this approach may produce imprecise numerical results or incur a high computational cost.

VIII. SIMULATION RESULTS

*VIII.1. Load Modeling ZIP loads to ZP loads*

Constant current loads are split into the constant impedance and the constant power loads using the formula $S_j^{app} = \left(\frac{S_j^0}{2|v_j^0|}\right)|v_j|^2 + \left(\frac{S_j^0}{2}\right)|v_j^0|$. Power flow studies are performed on all the systems in the MATPOWER [27] to identify voltage caused errors caused by the approximation of the constant current loads. The errors remain less than 1% for all the systems. In the simulations, loads are equally distributed as *Z*, *I*, and *P* loads, i.e., one third of loads are *Z* loads, *I* loads, and *P* loads each at all the *PQ* nodes. Figure 3 dipicts the typical examples from the selected systems.

If the voltage magnitudes do not deviate more than 10%, relative errors in the loads, as shown in Figure 1, remain below 1%. As a result, the resulting voltages are also less than 1%. The numerical error in approximating loads from constant current loads, as depicted in Figure 1, falls within a range similar to that of the demand forecast. The voltage difference is normally minor as long as demand forecast is close to the actual demand. In general, the resulting voltages are determined by how the load is distributed across the network. If the input is close to a certain value, the system's output will remain close to the expected output, as illustrated in Figure 3.

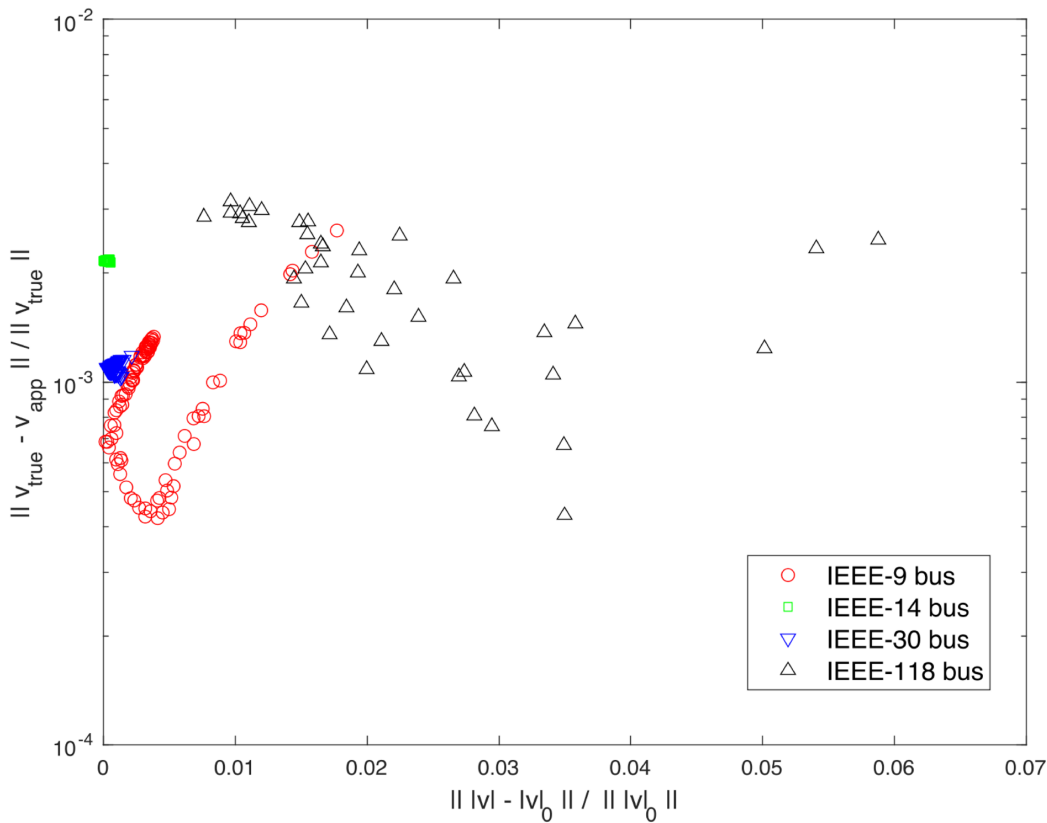

Figure 3. Errors in the voltages incurred by the modeling for the constant current loads in terms of the constant impedance loads and the constant power loads in power flow studies for IEEE-9, -14, -30, and -118 bus systems.

*VIII.2. Voltage Estimate Using Holomorphic Embedding*

Section IV discusses several circuit analysis methods for simplifying complex electrical networks. These methods represent an arbitrary two-port network is represented by these methods as a single voltage source in series with an impedance. By utilizing circuit analysis techniques, they enable determination of voltage and current behavior at various points in the network. The resulting equivalent circuit provides a mean to model the original circuit's behavior in terms of voltage, which proves valuable in reducing electrical circuit complexity and enhancing the understanding of their behavior. The networks considered for numerical simulation are the same as in Section VIII.1, with a typical example displayed in Figure 4. All tested networks perform similarly.

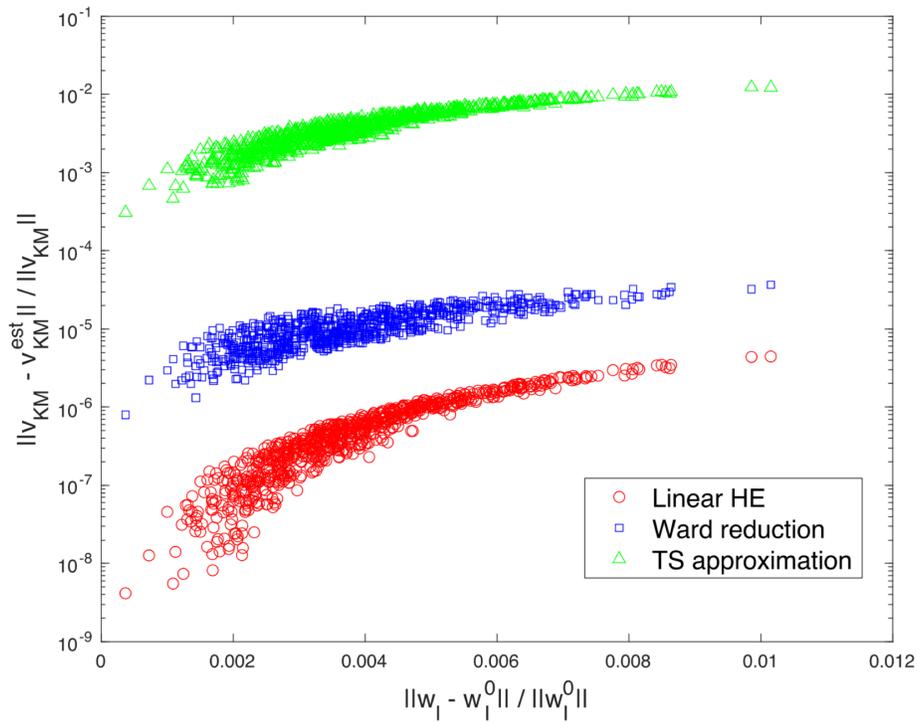

Figure 4. Errors in the voltage estimate at *KM* in the IEEE-30 bus system with respect to the deviation of the voltages at *Ibus* where the linearization are performed; 1) proposed linear holomorphic embedding approach (in red circle); 2) Ward reduction (in blue square); and 3) Truncated series approximation (in green triangle).

The proposed method clearly provides an accurate estimates of the voltages at *KM* while covering a wide range of voltages at *Ibus*. Considering the accuracy of the estimate and the formulation of all nodal voltages and branch flows are expressed in terms of voltages at *Ibus*, it is possible to construct an efficient, accurate, and highly adaptable "reduced" network that maps the power flows over the original network. However, the challenge of developing a network corresponding to the highly asymmetric nodal admittance matrix appearing in $I = Yv$ poses a potential barrier to creating a "reduced" network.

*VIII.2. Simulation Environments*

Simulations were conducted on various IEEE model systems (including the IEEE 4-, 9-, 14-, 30-, and 39-bus systems). However, for purpose of visual presentation and comparison with the results from [17], the simulation results for the IEEE 9-bus system are discussed in detail (Figure 5). The 9-bus system is expanded to 12 buses by adding 3 *Ibus* (Bus 10, 11, and 12) to represent the internal generators located at Bus 1, 2, and 3, respectively.

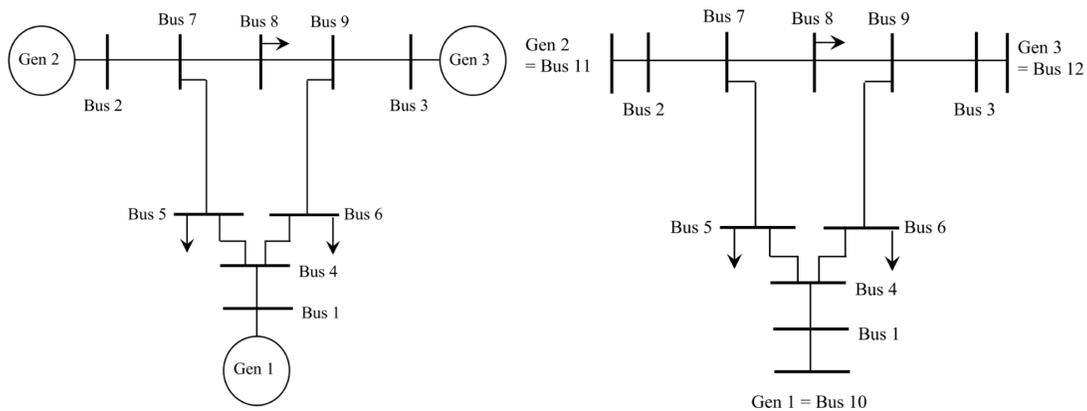

Figure 5. IEEE-9bus system shown in left. The arrows represent the system's loads, and they are all *ZIP* loads with equal distributions (1/3, 1/3, 1/3). Note that *Ibus* = {Gen 1, Gen 2, Gen 3}, *Kbus* = {Bus 1, Bus 2, Bus 3}, *Mbus* = {Bus 4, Bus 5, Bus 6, Bus 7, Bus 8, Bus 9}. An equivalent 12-bus system is built by modeling the generators as nodes, i.e., Gen 1, 2, and 3 correspond to Bus 10, 11, and 12, respectively (shown in right). Rotor angles are represented by the voltage angles at *Ibus*.

In one particular case, a line failure is assumed to occur near Bus 7 on the line connecting Bus 7 and Bus 8, which consist of two separate circuits. The fault is detected immediately after the disturbance, and the circuit breaker is opened to isolate the disturbance, indicating that the fault is cleared. Once the fault is cleared, power is restored to the remaining circuit. The stability is assessed using a coupled oscillator model [34], direct methods [2]–[7], [35], [36], numerical *TDS*, and the proposed method. Other outages, such as various load losses at Bus 8, are modelled and analyzed for their effects on system stability simulated and evaluated using the same methodologies.

When the voltage magnitudes at all nodes deviate from the pre-fault state along the on-fault trajectory of the line outage before the fault is cleared, the errors in load modeling, as well as the errors from the linearization in the holomorphic embedding, become considerable. In contrast to the errors discussed above, the inaccuracy induced by assuming that all loads are constant impedance loads is rather minor. Consequenty, utilizing the constant power load type during the on-fault trajectory is feasible. However, because the fault-clearing process is quick, precise *TDS* does not necessitate a high computation cost. When the loads are not entirely constant impedance loads, the on-fault trajectory is identified using *TDS* simulation in this study.

It is worth noting that the *HTMS* method converges within a time radius of around $10^{-2}$-$10^{-1}$ seconds, with 20-30 terms yielding accurate results when compared to *TDS* results, with a 5% tolerance between the two components. Each *HE* procedure at involves a matrix inversion, resulting in computation costs for *HTMS* being in the order of $\vartheta\left(\frac{m}{T_{HTMS}} Nb^3\right)$ [31], where *m* is the number of terms in the truncated McLaurin series and $T_{HTMS}$ is the time period of one *HTMS* evaluation. The time step $T_{TDS}$ for *TDS* in this study is $10^{-3}$, and *QPF* includes a matrix inversion, therefore the *TDS* computation costs are in the order of $\vartheta\left(\frac{1}{T_{TDS}} Nb^3\right)$ [31]. The costs of computing *HTMS* ($\pi_{HTMS}$) are comparable to those of computing *TDS* ($\pi_{TDS}$), $\frac{\pi_{HTMS}}{\pi_{TDS}} \propto \vartheta\left(\frac{T_{TDS}}{T_{HTMS}} m\right)$ [31].

All numerical calculations were performed on a Mac Pro equipped with two 2.93 GHz 6-core Xeon processors. The *P* load values are adjusted to match the *Z* load values in the post-fault trajectory presented in [17], i.e., $s_{Zload}^{j} = |v_j|^2 s_{base}^{j}$ and $s_{Pload}^{j} = s_{base}^{j}$, where $v_j$ denotes the voltage at Bus *j* at the beginning of the post-fault trajectory in [17]. Consequently, the stability assessment results in Ref [17] for the coupled oscillator and the direct method remain unchanged. Table 1 only shows the evaluation results for *TDS* and the proposed method.

Table 1. Summary of the simulation results

| Disturbance | Line Outage | | | | Loss of loads (10 %) | | | | Loss of loads (100 %) | | | |
|---|---|---|---|---|---|---|---|---|---|---|---|---|
| Model | TDS | | Proposed method | | TDS | | Proposed method | | TDS | | Proposed method | |
| Load model | Z loads | P loads | Z loads | P loads | Z loads | P loads | Z loads | P loads | Z loads | P loads | Z loads | P loads |
| Comp. time in sec | 3.37 | 5.44 | 1.91 | 3.71 | 3.80 | 6.33 | 0.17 | 1.15 | 2.62 | 4.35 | 0.55 | 3.77 |

*VIII.3. Line Outage Between Bus 5 and Bus 7 near Bus 7*

The line outage results will serve as an example. As shown in Table 1, the *TDS* approach is employed to calculate the on-fault trajectory, which increases the computation time of the proposed method,. In contrast to *Z* loads where the loads gradually adjust based on voltage magnitudes, the rotor motions are constrained to support constant power loads. Because of these constraints imposed on *P* loads, the real and imaginary components of the rotor motion (the time derivatives of the voltages) deviate significantly from zero. Voltages at *Ibus* (rotor) with *P* loads, as a result, cross the validity region boundaries more frequently than voltages at *Ibus* (rotor) with *Z* loads do to satisfy (17).

The validity region is a key element in the proposed method since it assures that the resulting solution obtained is sufficiently close to the time series provided by *TDS*. In the case of the line outage simulation, the validity region boundary $\left\| \left( 1 - \frac{1}{E_j}\sqrt{x_j^2 + y_j^2}; x_j\frac{dx_j}{dt} + y_j\frac{dy_j}{dt}; \frac{p_G}{p_G^0} - 1 \right) \right\| \cong \left\| x_j\frac{dx_j}{dt} + y_j\frac{dy_j}{dt} \right\| \leq \varepsilon$ is

reached at 0.1 seconds, as shown in 6. Figure 7 illustrates the real and imaginary components of the rotor motions as a result of position and velocity projections. Figure 7(C) compares the blue and red lines showing that when the value of $\varepsilon$ is set to 1% the approximation employed to simplify O1 holds well. The 1-norm is used to derive a useful expression for the validity boundary based on the inequality of norms [31], $\|\cdot\|_\infty \leq \|\cdot\|_2 \leq \|\cdot\|_1$.

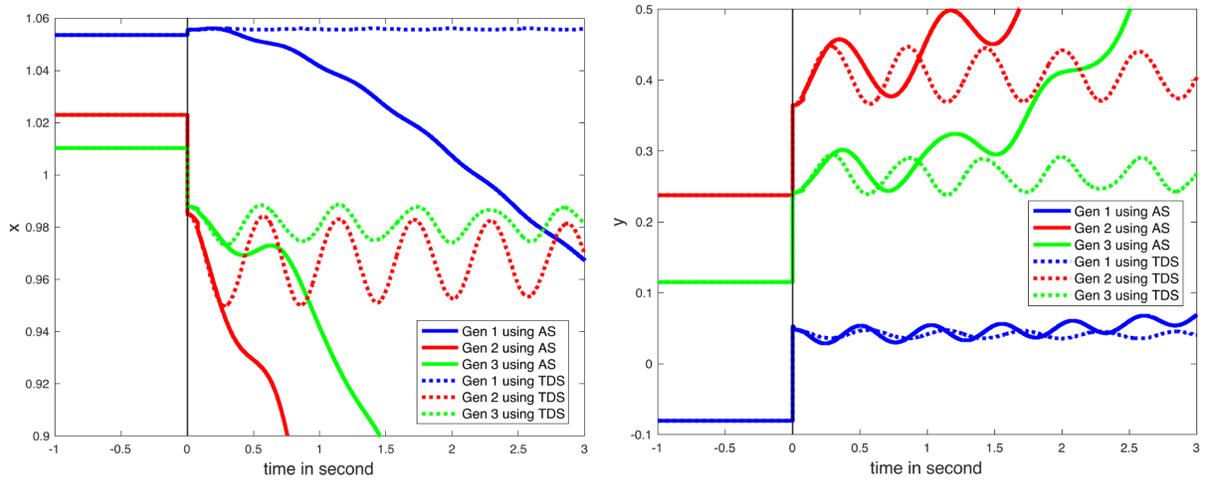

Figure 6. The behaviors after a line outage between Bus 5 and Bus 7 of (A) the real part of the rotor motion; (B) the imaginary part of the rotor motion.

The time where the time series intersects the validity boundary has an analytic form, and the consistent initialization yields $\rho_I, \omega_I,$ and $v_I$. Using the consistent initial values, the swing equations and corresponding solutions are identified for the next validity region, ensuring that the analytical solution and numerical computation remain sufficiently close to evaluate the stability over the entire range. The stability assessment is as follows. With a known analytical solution valid within a boundary, one can assess if the system is stable or becomes unstable within the boundary. If the system becomes unstable, the assessment determines that the system is unstable. Otherwise, the assessment is continued to the next boundary. If the assessment is not completed by the time frame that a control mechanism, such as governor control, is engaged, the system is at least operationally stable [17]. Only the value for $\beta$ in (18) should be modified to match the specified values for analytic solution. The computation costs associated with the update of $\beta$ are minimal because the evaluation of is a least square problem with a fixed matrix.

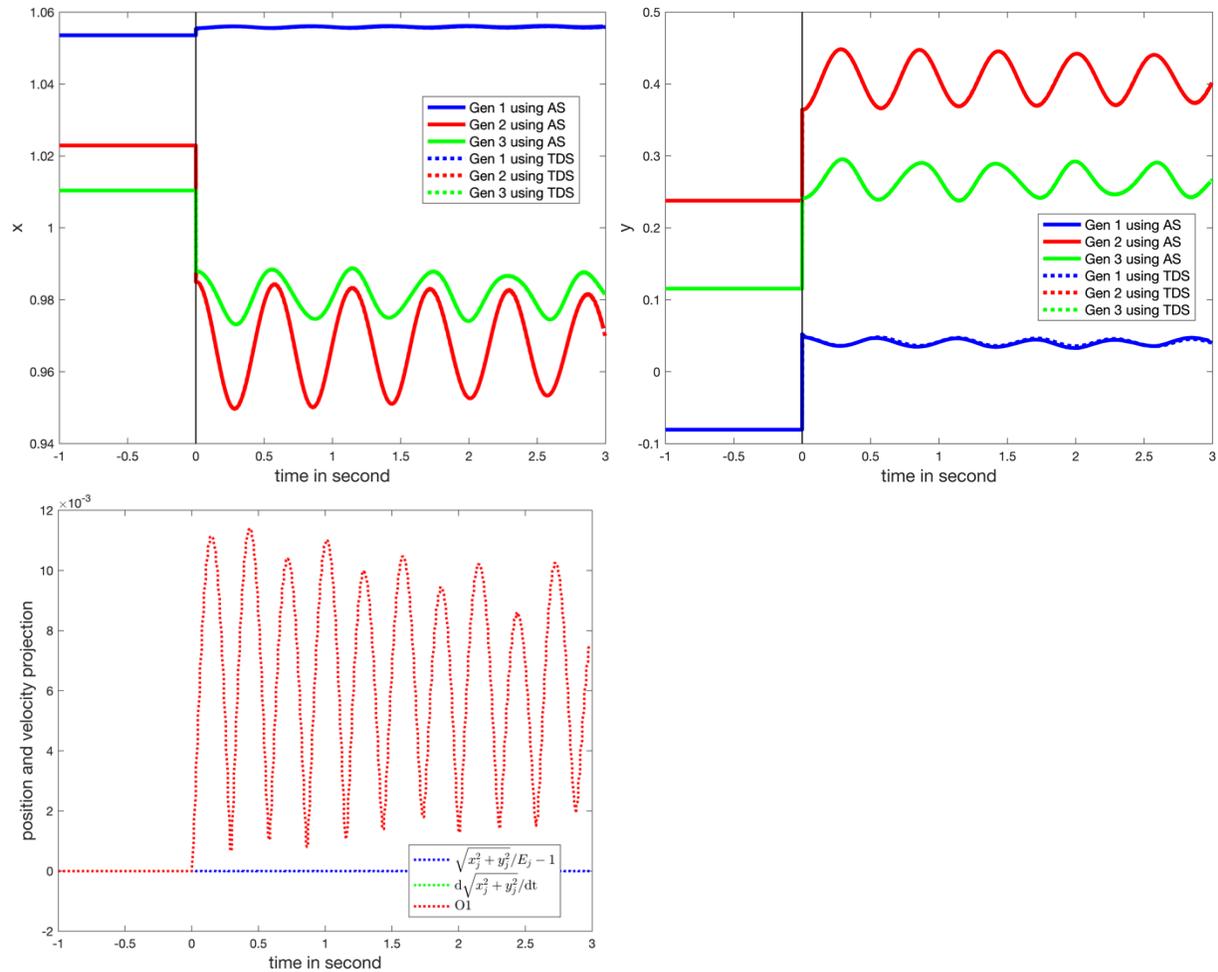

Figure 7. The behaviors with the consistent initializations after a line outage between Bus 5 and Bus 7 of (A) the real part of the rotor motion; (B) the imaginary part of the rotor motion; and (C) position projection (blue), velocity projection (green), and both (red).

It is important to note that the computation time in Table 1 includes all procedures until stability is determined, which entails approximately 10 times of consistent initialization problems until $t = 3$ seconds. This considerably increases the computation cost. The proposed approach, however, has a substantially reduced computation cost when compared to *TDS*. A greater value for $\varepsilon$ in (17) can be utilized to further improve computational efficiency at the expense of compromising analytical solution accuracy. Table 1 shows the computation times for *TDS* and the proposed method until stability is assessed ($t = 1$ second). When $\varepsilon$ is set to 10%, the initialization problem is only solved twice, while the stability assessment remains unaltered, resulting in a reduction in computation time from 3.71 seconds to 1.22 seconds.

*VIII.4. Voltage Linearization Using Holomorphic Embedding*

The *HE* linearization method relies on a specified reference angle, typically zero, at an arbitrarily chosen node (in this study, Bus 1). To achieve this, we rotate the voltages until the angle at the reference node is zero. Following the fault clearance ($t = 0$ sec), the rotor angles undergo step changes, as shown in Figure 8 (A) and (B). These changes result from the voltage angle shift at the reference angle node, not from alteration in the physical rotor angles. The rotor angles with respect to that of *COI* offsets the step change as illustrated in (20) (See Figure 8 (B)). *TDS* will continue until the system establishes a new equilibrium or becomes unstable. Each *TDS* step necessitates solving 2 *QPF* problems – one at the predictor stage and the other at the corrector stage, for a total of 6,000 *QPF* are solved during the 3-second at a frequency of $10^{-3}$ sec. As the proposed method offers an analytic expression for the variables over time during the post-fault period, it is possible to calculate errors. For example, the leading terms in the linearization error based on *HE* are of the order $\vartheta(NKM^2)$ [31] if $(Y_{tr}^{KMKM})^{-1}$ is stored.

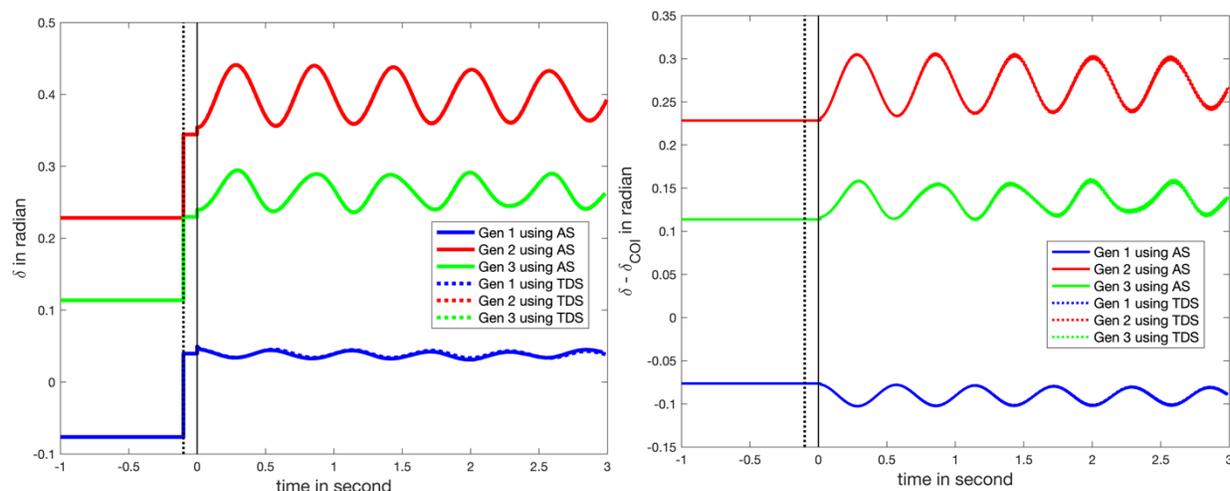

Figure 8. The rotor angles with the consistent initializations after a line outage between Bus 5 and Bus 7. The angles are referenced against the voltage angle (1) at Bus 1 or (2) at the Center of Inertia (COI).

Figure 9 illustrates nodal voltages in the IEEE 9-bus system, where the imaginary components of the voltages at selected nodes (Bus 1, 4, 5, and 6) are non-positive. This observation implies that the behavior of the nodal voltages may not be uniform across the system, indicating that there may be an underlying structure that warrants further investigation. To investigate this further, the sensitivity matrix $H_{KM}$, as given in (10), is examined. Let $H_{KM}^j$ represent the $H_{KM}$ sensitivities corresponding to Bus $j$. The angle between two subspaces formed by $H_{KM}^i$ and $H_{KM}^j$ is defined as $\theta_{ij}$. Two spaces are considered nearly linearly dependent when the angle between them is small [37],[38] (See Table 2 for the values for $\theta_{ij}$). Partitioning the buses based on the behavior exhibited in Figure 9 may be difficult because $\theta_{ij}$ represents the angle between two subspaces spanning 2$NI$ dimensions. To cluster the buses, the *k*-means clustering algorithm [39],[40] is employed, with the $\theta_{ij}$ values. According to *k*-means clustering results, Buses 1, 4, 5, and 6 constitute one group, while the remaining buses form another. Figure 5 displays the locations of these groups in the IEEE 9-bus system, with the lower half of the system forming one group. Following a line outage, the path between the two groups becomes less efficient for power transfer, leading to a degree of system partitioning. This partitioning impact is reflected in the sensitivity matrix $H_{KM}$, which can provide valuable insights for analyzing the system's stability and vulnerability during disturbances. Further studies can focus on understanding and mitigating the impacts of these partitioning effects on power system performance and reliability.

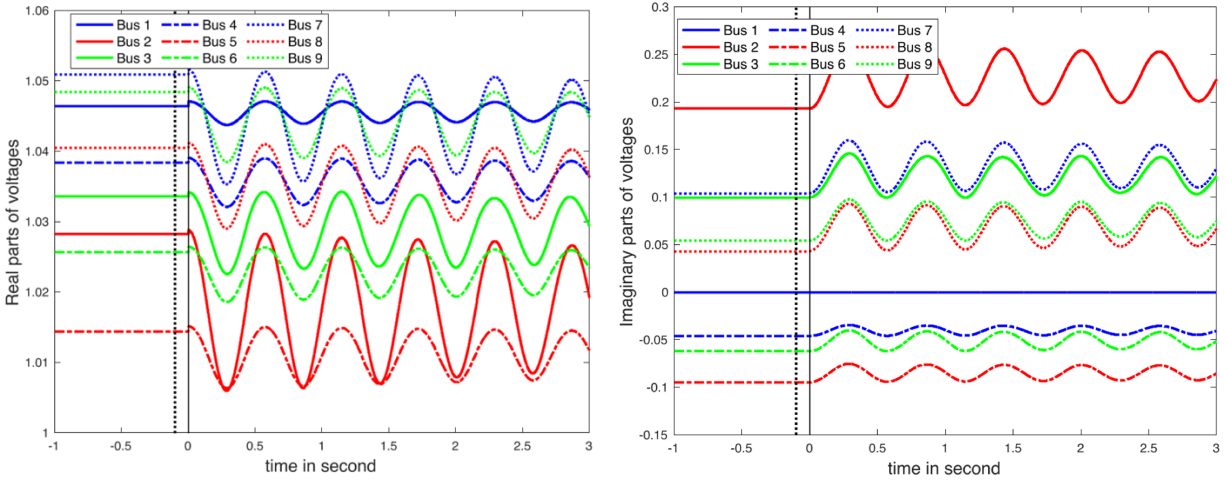

Figure 9. The timeseries of nodal voltages after a line outage between Bus 5 and Bus 7: (A) the real parts of the nodal voltages; and (B) the imaginary parts of the nodal voltages until $t$ = 3 seconds

.

Table 2. The angles in radian between two subspaces, $\theta_{ij}$

| $\theta_{ij}$ | Bus 1 | Bus 2 | Bus 3 | Bus 4 | Bus 5 | Bus 6 | Bus 7 | Bus 8 | Bus 9 |
|---|---|---|---|---|---|---|---|---|---|
| Bus 1 | - | 1.126 | 0.988 | 0.171 | 0.303 | 0.368 | 0.931 | 0.863 | 0.818 |
| Bus 2 | 1.12 | - | 1.079 | 1.049 | 0.877 | 0.928 | 0.243 | 0.437 | 0.866 |
| Bus 3 | 0.988 | 1.079 | - | 0.930 | 0.876 | 0.675 | 0.915 | 0.653 | 0.239 |
| Bus 4 | 0.171 | 1.049 | 0.930 | - | 0.253 | 0.312 | 0.880 | 0.794 | 0.771 |
| Bus 5 | 0.303 | 0.877 | 0.876 | 0.253 | - | 0.403 | 0.692 | 0.622 | 0.684 |
| Bus 6 | 0.368 | 0.928 | 0.675 | 0.311 | 0.403 | - | 0.738 | 0.583 | 0.499 |
| Bus 7 | 0.931 | 0.243 | 0.915 | 0.880 | 0.692 | 0.738 | - | 0.263 | 0.701 |
| Bus 8 | 0.863 | 0.437 | 0.653 | 0.794 | 0.622 | 0.583 | 0.263 | - | 0.438 |
| Bus 9 | 0.818 | 0.866 | 0.239 | 0.771 | 0.684 | 0.499 | 0.701 | 0.438 | - |

*VIII.5. The Effects of Load Characteristics on Power System Dynamics*

Figure 10 (A) displays the characteristics of various load types, including constant impedance loads (*Z*), constant current loads (*I*), and constant power loads (*P*), using the proposed method, which yields results that are visually indistinguishable from those obtained using the TDS method. This implies that in this study,

the two major error sources of the proposed analytical solution to swing equations - *I* load modeling in terms of *Z* and *P* loads, and the condition in (17) - are minimal.

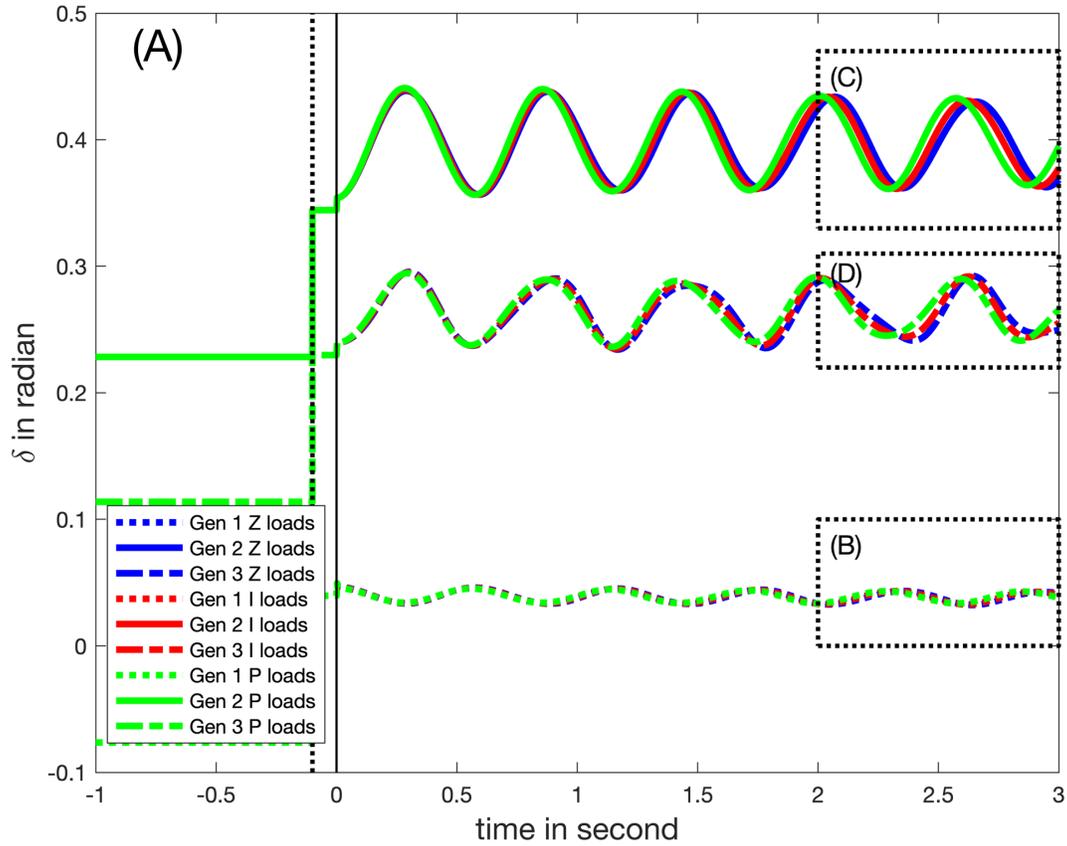

Figure 10. The timeseries of nodal voltages after a line outage between Bus 5 and Bus 7: (A) the rotor angles (the voltage phasor at *Ibus*); (B) the rotor angles at Generator 1 at $t \in [2, 3]$; (C) the rotor angles at Generator 2 at $t \in [2, 3]$; (D) the rotor angles at Generator 3 at $t \in [2, 3]$. Blue, red, and green lines correspond to the load types, constant impedance (*Z*), constant current (*I*), and constant power (*P*) loads, respectively.

Figure 10 (A) initially displays load characteristics before 0.5 seconds, which become evident after 2 seconds (for further details, see Figures 10 (B), (C), and (D)). According to Section III.3, higher nodal voltages at a node necessitate a larger power for P loads. This is reflected in $Lw_I$ in (15), resulting in a larger *L* matrix.

The analysis of the eigenvalues and coefficients $\beta_k^* \Psi_k$ from (18) reveals that only two coefficients associated with the generator at Bus 3 have significant values (see second row of Table 3). Their corresponding eigenvalues are listed in the third and the fourth rows of Table 3. The first set of eigenvalues have larger imaginary components, suggesting that the corresponding functions $u$ in (18) oscillate more quickly than those related to the second set. Therefore, the smaller coefficients of the second set of imaginary eigenvalues cause the modes' linear sum to oscillate even faster. This explains why the periods of the rotor dynamics decreases with $I$ and $P$ loads, potentially due to the damping effect of the constant power loads on the system.

Table 3. The two largest coefficients and the corresponding eigenvalues associated with Generator 3

|  | Z loads | | I loads | | P loads | |
| --- | --- | --- | --- | --- | --- | --- |
| Coefficients, $\beta_k^* \Psi_k$ | -0.0255 | -0.0118 | -0.0259 | -0,0079 | -0.0264 | 0.0001 |
| Eigenvalues (real part) | -0.0947 | -0.1013 | -0.0945 | -0.1013 | -0.0946 | -0.1013 |
| Eigenvalues (imaginary part) | 7.6788 | -0.3682 | 7.7554 | -0.3720 | 7.8606 | -0.3791 |

To simplify our physical explanation, let's disregard $l$ in (15) and suppose that all coefficients in (15) are scalar. With this adjustment, (15) then represents a spring-mass-dashpot system, where a spring either push or pull the mass illustrated in Figure 11. The equation of the motion for this system is: $m \frac{d^2}{dt^2} x + c \frac{d}{dt} x + kx = 0$ [33] where $m$ denotes the mass, $c$ is the dashpot constant, and $k$ represents the spring constant. The solution to this equation is: $x(t) = x_1 exp(\lambda_1 t) + x_2 exp(\lambda_2 t)$ where $\lambda_1 = \frac{1}{2m}\left(-c + \sqrt{c^2 - 4mk}\right)$ and $\lambda_2 = \frac{1}{2m}\left(-c - \sqrt{c^2 - 4mk}\right)$. The values for $x_1$ and $x_2$ are selected to satisfy the system's initial and boundary conditions. When $c^2$ is less than $4mk$, the system oscillates underdamped.

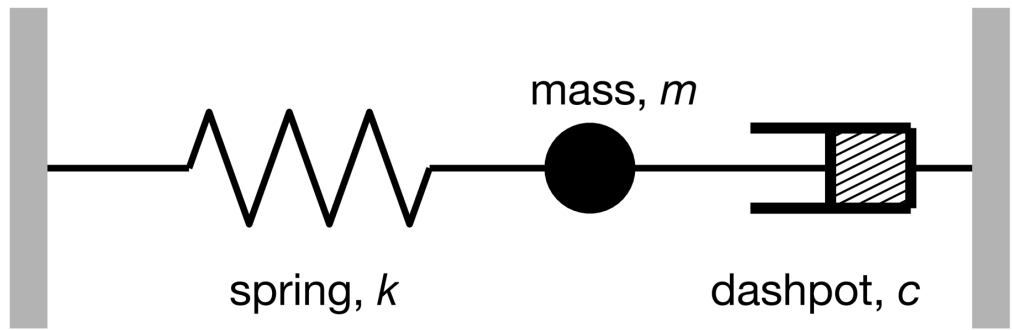

Figure 11. A schematic diagram for a spring-mass-dash system where *k* is the spring constant; *m* is the mass; and *c* is the dashpot constant.

Figure 12 demonstrates different damping behaviors based on various k values, assuming $x_1 = x_2 = 1$, $m = 1$, and $c = 1.7$. In case where $k = 1$ (represented by the blue line), the system is overdamped since $c^2 > 4mk$, and the mass shifts to the next stable point. The red and green lines display changes in the mass location for lower and higher k values respectively. Importantly, Figure 12 shows that as *k* value increases, the oscillation period decreases.

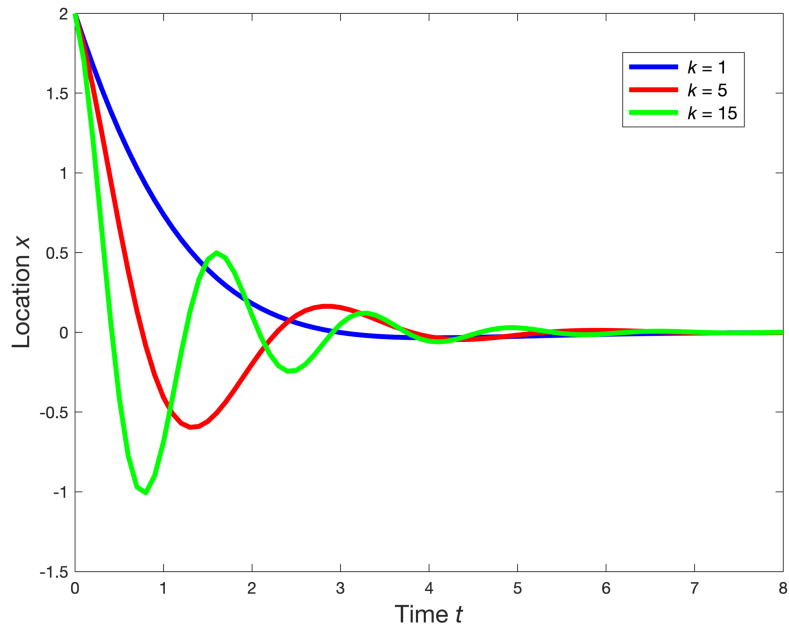

Figure 12. The motions of the mass according to various values for the spring constant *k*; *k* = 1 (Blue); *k* = 5 (Red); and *k* = 15 (Green).

Figure 10 shows that the oscillation periods $\tau$ change with the types of loads, $\tau_{Zload} > \tau_{Iload} > \tau_{Pload}$. This suggests that in the spring-mass-dashpot system, the loads function similarly to a spring, and constant power loads the strongest restoring force.

IX. CONCLUSIONS

The authors of this paper propose a novel approach for evaluating the transient stability of power grids by deriving an analytic solution to the swing equation that can accommodate all ZIP models while complying to Kirchhoff's laws. The swing equation is transformed into a differential algebraic equation, and that the ZIP load models are analytically integrated using the holomorphic embedding method. As demonstrated in the paper's example, the proposed solution enables an accurate and efficient assessment of stability.

The proposed assessment tool is fundamentally different from traditional methods such as the coupled oscillator model and the direct method because it avoids making physically inadmissible assumptions. Additionally, it distinguishes itself from time-domain simulation methods by offering relatively low computation costs. These advantages make the proposed tool an attractive option for real-time transient stability assessment in power systems, which is crucial for maintaining grid reliability and security.

The potential applications of this novel method include real-time monitoring and control of power systems, as well as grid planning and optimization. By providing a more accurate and computationally efficient means of assessing transient stability, the proposed approach can help grid operators make better-informed decisions and respond more quickly to potential threats, ultimately enhancing the overall performance and resilience of power systems.

**Author Contributions**

Conceptualization: HyungSeon Oh.

Data curation: HyungSeon Oh.

Formal analysis: HyungSeon Oh.

Investigation: HyungSeon Oh.

Methodology: HyungSeon Oh.

Software: HyungSeon Oh.

Validation: HyungSeon Oh.

Writing – original draft: HyungSeon Oh.

Writing – review & editing: HyungSeon Oh.